# LOWER ESTIMATES OF TRANSITION DENSITIES AND BOUNDS ON EXPONENTIAL ERGODICITY FOR STOCHASTIC PDE'S

By B. Goldys and B. Maslowski[1]

*University of New South Wales and Academy of Sciences of Czech Republic*

A formula for the transition density of a Markov process defined by an infinite-dimensional stochastic equation is given in terms of the Ornstein–Uhlenbeck bridge and a useful lower estimate on the density is provided. As a consequence, uniform exponential ergodicity and $V$-ergodicity are proved for a large class of equations. We also provide computable bounds on the convergence rates and the spectral gap for the Markov semigroups defined by the equations. The bounds turn out to be uniform with respect to a large family of nonlinear drift coefficients. Examples of finite-dimensional stochastic equations and semilinear parabolic equations are given.

**1. Introduction.** The aim of this paper is to study the ergodic properties of solutions to a semilinear stochastic equation

$$\begin{aligned} dX^x &= (AX^x + F(X^x))\,dt + \sqrt{Q}\,dW, \\ X_0^x &= x \in E, \end{aligned} \tag{1.1}$$

considered in a separable Banach space $E$, where $W$ is a cylindrical Wiener process on a Hilbert space $H$ such that $E \subset H$. Under the assumptions listed below (see Section 2), this equation has a unique Markov solution $(X_t^x)$ with a unique invariant measure $\mu^*$.

Ergodic properties of solutions to infinite-dimensional stochastic differential equations have been extensively studied in recent years. The key problems in this field are the existence and uniqueness of invariant measure and the rate of convergence of the time $t$ distribution of the process to the invariant measure. In the case of $\dim E < \infty$ these questions have

Received February 2004; revised June 2005.
[1]Supported in part by the ARC Discovery Grant DP0346406, the UNSW Faculty Research Grant PS05345 and the GAČR Grants 201/01/1197 and 201/04/0750.

been studied for a long time and the ergodic theory of finite-dimensional diffusion processes is relatively well developed, see, for example, a classical monograph [21]. In this paper we study the ergodic properties of a class of ordinary and partial stochastic differential equations that includes stochastic reaction–diffusion equations in bounded domains. First results on the existence and uniqueness of invariant measures for stochastic reaction-diffusion equations were obtained in [12, 26, 42], see also [5], the monographs [6, 10] and references therein. The rate of convergence to the invariant measure in infinite dimensions became a subject of interest much later and still is not well understood. Jacquot and Royer [24] proved exponential ergodicity for a semilinear parabolic equation with bounded nonlinear drift, Shardlow [39] applied the theory of Meyn and Tweedie to obtain $V$-uniform ergodicity for some semilinear equations in Hilbert spaces. Hairer in [14] proved, under different sets of conditions, uniform exponential ergodicity for equations with drifts growing faster than linearly. Exponential convergence to equilibrium in a norm intermediate between the total variation metric and the Wasserstein metric has been obtained in [31] for the stochastic Navier–Stokes equation.

A closely related problem of asymptotic behavior of the Markov semigroup $P_t\phi(x) = \mathbb{E}\phi(X_t^x)$ attracted much attention due to its importance in Mathematical Physics. In particular, exponential convergence of the semigroup in the spaces $L^p(E, \mu^*)$, $p \in [1, \infty)$, and related questions of the existence of the spectral gap and logarithmic Sobolev inequality have been studied by numerous authors, see [1, 2, 3, 7, 8, 11, 18, 22, 43].

The aim of the present paper is to prove $V$-uniform (exponential) ergodicity with $V(x) = |x|_E + 1$ and, if the drift grows faster than linearly, uniform exponential ergodicity, for equation (1.1). Our method allows us to find exact bounds on convergence (i.e., to give explicit estimates for the rate of exponential convergence in the total variation norm or $V$-variation norm). In this respect, our results seem to be new even for finite-dimensional SDE's (which is also due to our method to estimate the transition density that, to the best of our knowledge, has not been used in finite dimensions so far). If the Markov semigroup $(P_t)$ is symmetric, we obtain explicit lower estimates for the spectral gap in $L^2(E, \mu^*)$. Stronger results are obtained in the case of a drift growing faster than linearly: for a symmetric Markov semigroup, we show uniform estimates on the spectral gap in the spaces $L^p(E, \mu^*)$ for all $p \in [1, \infty)$ and in the nonsymmetric case, our estimates remain valid for $p > 1$, in particular, in $L^2(E, \mu^*)$.

Unlike in the aforementioned papers, in the present paper a lower bound measure and a suitable small set for a skeleton process are found explicitly in terms of the lower estimates of transition densities and the constants in an ultimate boundedness condition (or, in particular, a suitable Lyapunov function). This enables us to apply earlier results on computable bounds for Markov chains, which are expressed in terms of lower bound measures, corresponding small sets and constants from the Lyapunov–Foster geometric drift condition [33]. The bounds turn out to be uniform with respect to a large family of drift coefficients, which is important for proving continuous dependence of invariant measures on parameter (cf. Section 8). We also believe that this uniformity is an important tool for studying the Hamilton–Jacobi–Bellman equation for the ergodic control problem. On the other hand, the method employed here has its limitations. Our method strongly relies on the Girsanov theorem and therefore, we need an assumption that $F$ maps the whole state space into the range of $\sqrt{Q}$. Therefore, any extension to other types of equations (like stochastic Burgers or Navier–Stokes equations) would be difficult. Note, however, that in the two recent authors' papers [17] and [16] $V$-uniform ergodicity and spectral gap type results have been proved for stochastic Burgers, 2D Navier–Stokes and more general reaction–diffusion equations. Nonetheless, in these papers a different method is used that allows us neither to find explicit bounds on the convergence constants nor to show the uniformity of convergence with respect to coefficients.

An important tool for our proofs is a formula for the transition densities that is derived in this paper. We use this formula to establish suitable lower estimates on the densities which we believe are of independent interest. They are obtained by means of the Girsanov theorem and the so called Ornstein–Uhlenbeck bridge (or pinned Ornstein–Uhlenbeck process). Let us explain the main idea of this approach.

Let $(Z_t^x)$ be an Ornstein–Uhlenbeck process on a separable Hilbert space $H$. By this, we mean that $(Z_t^x)$ is a solution to a linear stochastic evolution equation

$$
\begin{aligned}
dZ_t^x &= AZ_t^x\,dt + \sqrt{Q}\,dW_t, \\
Z_0^x &= x \in H.
\end{aligned}
\tag{1.2}
$$

The Ornstein–Uhlenbeck bridge $(\widehat{Z}_t^{x,y})$ associated to the Ornstein–Uhlenbeck process $(Z_t^x)$ is informally defined by the formula

$$\mathbb{P}(Z_t^x \in B | Z_1^x = y) = \mathbb{P}(\widehat{Z}_t^{x,y} \in B), \qquad t < 1,$$

where $x, y \in H$ and $B \subset H$ is a Borel set. The importance of various types of bridge processes for the study of transition densities of finite dimensional diffusions is well recognised, see, for example, [23]. In infinite-dimensional framework this concept was developed in [41] in order to study the regularity of transition semigroups of diffusions on Hilbert spaces. In [28] and [29] an Ornstein–Uhlenbeck bridge is introduced in order to obtain lower estimates on the transition kernel of some semilinear stochastic evolution equations. The basic idea is as follows. Using the equivalence of measures corresponding to $X_t^x$ and $Z_t^x$ and the Girsanov formula, we can write the transition density of the process $X_t^x$ in the form

$$d(T, x, y) = \mathbb{E}(\Phi(Z_\cdot^x)|Z_T^x = y),$$

where $\Phi$ is a measurable functional defined on trajectories of the Ornstein–Uhlenbeck process. This form of the density is not suitable for the uniform estimates that are needed. Therefore, the conditional expectation is transformed into a usual expectation with respect to the measure of the OU bridge $(\widehat{Z}_t^{x,y})$ considered for $t \in [0, T]$:

$$d(T, x, y) = \mathbb{E}\Phi(\widehat{Z}_\cdot^{x,y}),$$

which enables us to find the uniform lower estimates. Let us note here a technical difficulty caused by the fact that we can define the OU bridge for $y$ in a certain Borel subspace of measure one only, but this turns out to be sufficient for our needs.

Precise formulations and hypothesis are given in the following Section 2. In Sections 3 and 4 the properties of the OU bridge, which are needed in the sequel, are established (some auxiliary results are deferred to the Appendix). The formula for transition densities is found and the lower estimates are given in Section 5. These results are applied in Section 6 to establish our main results, $V$-uniform ergodicity and uniform exponential ergodicity, respectively, and find the computable bounds on respective constants. In Section 7 the corollaries on the $L^p(E, \mu^*)$ exponential convergence and the spectral gap are stated. Section 8 is devoted to some extensions and applications (continuous dependence of invariant measures on a parameter). Examples (finite-dimensional nonlinear stochastic oscillator and stochastic parabolic equations) are presented in Section 9.

**2. Assumptions and notation.** Let $H = (H, |\cdot|)$ be a real separable Hilbert space and let $E = (E, |\cdot|_E)$ be a separable Banach space densely embedded into $H$. In this paper we will study a stochastic semilinear equation

$$dX_t = (AX_t + F(X_t))\,dt + \sqrt{Q}\,dW_t,$$
(2.1)
$$X_0 = x \in E,$$

where $(W_t)$ is a standard cylindrical Wiener process on $H$ defined on a stochastic basis $(\Omega, \mathcal{F}, (\mathcal{F}_t), \mathbb{P})$ satisfying the usual conditions, $A$ denotes a linear operator on $H$ generating a strongly continuous semigroup $(S_t)$ on $H$ and $F$ is a nonlinear mapping $E \to E$. The first assumption assures the existence of an $H$-valued and strong Feller solution to the linear version of (2.1), when $F = 0$; in this case we consider the linear equation

$$dZ_t^x = AZ_t^x\,dt + \sqrt{Q}\,dW_t,$$
(2.2)
$$Z_0^x = x.$$

The solution to equation (2.2) is given by formula

$$Z_t^x = S_t x + \int_0^t S_{t-s} \sqrt{Q}\, dW_s, \qquad t \geq 0. \tag{2.3}$$

HYPOTHESIS 2.1. The operator $Q \geq 0$ is bounded and symmetric. For each $t > 0$, a bounded operator,

$$Q_t = \int_0^t S_s Q S_s^* \, ds,$$

is of trace class. Moreover,

$$\operatorname{im}(S_t) \subset \operatorname{im}(Q_t^{1/2}), \qquad t > 0. \tag{2.4}$$

If Hypothesis 2.1 holds, then

$$\overline{\operatorname{im}(Q_t)} = H, \qquad t > 0. \tag{2.5}$$

It is well known (cf. [9]) that (2.4) is equivalent to the strong Feller property of the process $(Z_t^x)$. Moreover, (2.4) yields

$$\operatorname{im}(Q_t^{1/2}) = \operatorname{im}(Q_1^{1/2}), \qquad t > 0. \tag{2.6}$$

The next hypothesis assures that the Ornstein–Uhlenbeck process $Z^x$ defined by equation (2.3) takes values in the Banach space $E$ and is continuous in $E$.

HYPOTHESIS 2.2. (a) The part $\tilde{A}$ of $A$ in the space $E$,

$$\tilde{A} = A|\operatorname{dom}(\tilde{A}), \qquad \operatorname{dom}(\tilde{A}) = \{y \in \operatorname{dom}(A) \cap E : Ay \in E\},$$

generates a $C_0$-semigroup on $E$, which is again denoted by $(S_t)$.
(b) The process $Z^0$ is $\mathbb{P}$-a.s. $E$-valued and $E$-continuous.

We further assume the following:

HYPOTHESIS 2.3.

$$\int_0^1 \|Q_t^{-1/2} S_t Q^{1/2}\|_{\mathrm{HS}}\, dt < \infty, \tag{2.7}$$

where $\|T\|_{\mathrm{HS}}$ stands for the Hilbert–Schmidt norm of the operator $T$.

Assumption (2.7) is not standard. It is needed to obtain a formula for the transition density (cf. Theorem 5.2). We will show that it is satisfied in many important cases (cf. Lemma 3.3, Remark 3.4 and Section 9).

In this paper we consider mild pathwise continuous solutions of (2.1). A process $X$ defined on a filtered probability space $(\Omega, \mathcal{F}, (\mathcal{F}_t), \mathbb{P})$ is a solution to equation (2.1) on an interval $[0,T]$ if $\mathbb{P}(X. \in C(0,T:E)) = 1$ and

$$(2.8) \quad X_t = S_t x + \int_0^t S_{t-r} F(X_r)\, dr + \int_0^t S_{t-r} \sqrt{Q}\, dW_r, \qquad t \in [0,T], \mathbb{P}\text{-a.s.}$$

Now we will formulate assumptions involving the nonlinear term $F$ in equation (2.1).

HYPOTHESIS 2.4. (a) The mapping $F: E \to E$ is Lipschitz continuous on bounded sets of $E$. For eaxh $x \in E$, there exists a unique mild solution $X$ to equation (2.1). Moreover, $X$ is a Markov process in $E$.

(b) $\text{im}(F) \subset \text{im}(Q^{1/2})$ and there exists a continuous function $G: E \to H$ such that $Q^{1/2} G = F$ and for some constants $K, m > 0$,

$$(2.9) \qquad |G(x)| \leq K(1 + |x|_E^m), \qquad x \in E.$$

REMARK 2.5. The assumption of local Lipschitz continuity of the mapping $F$ is not necessary for our main results. It may be replaced by the existence and uniqueness conditions for equation (2.1) and approximating equations. Similarly, the mapping $G$ need not be continuous. Only measurability and the polynomial bound (2.9), are needed but this would make some proofs technically more complicated.

Hypotheses 2.1–2.4 are standing assumptions of the paper and the results will be enunciated without recalling them again. Obviously, the hypotheses are used selectively (e.g., Hypothesis 2.4 is not needed for results on the Ornstein–Uhlenbeck bridge).

Denote by $\mathcal{B}$, $\mathcal{P}$ and $b\mathcal{B}$, the Borel $\sigma$-algebra of $E$, the space of probability measures on $E$ and the space of bounded Borel functions on $E$, respectively. Furthermore,

$$P_t \varphi(x) := \mathbb{E}_x \varphi(X_t), \qquad \phi \in b\mathcal{B}, x \in E, t \geq 0,$$

and

$$P(t, x, \Gamma) := P_t \mathbf{1}_\Gamma(x), \qquad x \in E, \Gamma \in \mathcal{B}, t \geq 0.$$

Let $(P_t^*)$ denote the adjoint Markov semigroup, that is,

$$(2.10) \qquad P_t^* \nu(\Gamma) := \int_\Gamma P(t, x, \Gamma) \nu(dx), \qquad t \geq 0, \nu \in \mathcal{P}, \Gamma \in \mathcal{B}.$$

An invariant measure $\mu^* \in \mathcal{P}$ is defined as a stationary point of the semigroup $(P_t^*)$, that is, $P_t^* \mu^* = \mu^*$ for each $t \geq 0$. Obviously, $P_t^* \nu$ is interpreted as the probability distribution of $X_t$ if $X_0$ has the initial distribution is $\nu$.

In our main theorems on $V$-uniform ergodicity, exponential ergodicity and spectral gap the solution to equation (2.1) is supposed to be ultimately bounded. In order to illustrate which systems are covered, it may be useful to formulate a growth condition on the nonlinear term $F$ which will be selectively used in some statements below (though it is not needed in our general theorems). By $\langle \cdot, \cdot \rangle_{E,E^*}$, we denote the duality between $E$ and $E^*$ and by $\partial |\cdot|_E$, the subdifferential of the norm $|\cdot|_E$. Suppose that there exist $k_1, k_2, k_3 > 0$, and $s > 0$ such that, for $x \in \operatorname{dom}(\tilde{A})$ and $x^* \in \partial |x|_E$, we have

$$(2.11) \quad \langle \tilde{A}x + F(x+y), x^* \rangle_{E,E^*} \leq -k_1 |x|_E + k_2 |y|_E^s + k_3, \qquad y \in E.$$

For example, if the mapping $F: E \to E$ is Lipschitz continuous on bounded sets in $E$ and Hypotheses 2.1 and 2.2 are satisfied, then the above condition implies existence of a unique mild solution to the equation (2.1) [i.e., Hypothesis 2.4(a)]. If, moreover, the moments of the Ornstein–Uhlenbeck process $Z^0$ are bounded on $[0, \infty)$ [condition (6.1) below], then there exists an invariant measure for the corresponding Markov process.

**3. Some properties of the Ornstein–Uhlenbeck process.** We will use the notation $\mu_t^x$ for the probability distribution of $Z_t^x$ and $\mu_t$ if $x = 0$. Obviously, $\mu_t^x$ is a Gaussian measure $N(S_t x, Q_t)$. For simplicity of notation, we set $Z_s := Z_s^0$, $s \geq 0$. It is easy to check that, for $s \leq t$

$$(3.1) \qquad \mathbb{E}\langle Z_s, h \rangle \langle Z_t, k \rangle = \langle S_{t-s} Q_s h, k \rangle.$$

LEMMA 3.1. *The operator $V_t = Q_1^{-1/2} S_{1-t} Q_t^{1/2}$ is bounded on $H$ and*

$$(3.2) \qquad \|V_t\| < 1, \qquad t \in (0, 1].$$

*Moreover,*

$$(3.3) \qquad \lim_{t \to 1} V_t^* x = \lim_{t \to 1} V_t x = x, \qquad x \in H.$$

PROOF. Estimate (3.2) was proved in [35]. It follows from (3.2) and a simple identity

$$Q_1 = Q_{1-t} + S_{1-t} Q_t S_{1-t}^*,$$

that

$$(3.4) \qquad Q_{1-t} = Q_1^{1/2} (I - V_t V_t^*) Q_1^{1/2}.$$

To prove (3.3), we will show first that

$$(3.5) \qquad \lim_{t \to 0} \langle V_t x, y \rangle = \langle x, y \rangle, \qquad x, y \in H.$$

Indeed, for $y \in \mathrm{im}(Q_1^{-1/2})$, we have

$$\lim_{t \to 1} \langle V_t x, y \rangle = \lim_{t \to 1} \langle S_{1-t} Q_t^{1/2} x, Q_1^{-1/2} y \rangle = \langle x, y \rangle.$$

For arbitrary $y \in H$, we may find a sequence $(y_n) \subset \mathrm{im}(Q_1^{-1/2})$, such that $y_n \to y$ in $H$ and then (3.2) yields

$$\langle V_t x, y - y_n \rangle \to 0,$$

uniformly in $t \leq 1$, and (3.5) follows. Next, (3.4) yields

$$\langle Q_{1-t} x, x \rangle = \langle (I - V_t V_t^*) Q_1^{1/2} x, Q_1^{1/2} x \rangle, \qquad x \in H.$$

It follows that, for each $y \in \mathrm{im}(Q_1^{1/2})$, we have

$$\lim_{t \to 1}(|y|^2 - |V_t^* y|^2) = 0,$$

and since $\|V_t^*\| < 1$ for all $t$, we find that

$$\lim_{t \to 1} |V_t^* y| = |y|, \qquad y \in H.$$

Now, invoking (3.5), we obtain the first part of (3.3). It is enough to prove the second part of (3.3) for $x$ such that $|x| = 1$. In this case (3.5) implies $\langle V_t x, x \rangle \to 1$, and thereby, invoking (3.2),

$$1 = \liminf_{t \to 1} \langle V_t x, x \rangle \leq \liminf_{t \to 1} |V_t x| \leq \limsup_{t \to 1} |V_t x| \leq 1.$$

Therefore, $V_t x \to 1$ as $t \to 1$. Now, taking into account (3.5), we obtain the second part of (3.3). $\square$

Clearly, $V_t^* = \overline{Q_t^{1/2} S_{1-t}^* Q_1^{-1/2}}$ and the operator

(3.6) $$K_t := Q_t^{1/2} V_t^*$$

is of Hilbert–Schmidt type on $H$. Then the operator

$$H \ni x \to \mathcal{K}x(t) := K_t x \in L^2(0, 1; H)$$

is bounded.

Let $\mu$ denote the probability distribution of the process $\{Z_t, t \in [0,1]\}$ concentrated on $L^2(0,1;H)$ and let $\mathcal{L}: L^2(0,1;H) \to C(0,1;H)$ be defined by the formula

$$\mathcal{L}u(t) = \int_0^t S_{t-s} Q^{1/2} u(s) \, ds.$$

The space $\mathrm{im}(\mathcal{L})$ endowed with the norm

$$\|\phi\| = \inf\{|u| : u \in L^2(0,1;H), \mathcal{L}u = \phi\}$$

may be identified with reproducing kernel Hilbert space of the measure $\mu$, see [9]. For any $t \in [0, 1)$, we define an unbounded $H$-valued operator

$$B_t x = Q^{1/2} S^*_{1-t} Q_1^{-1/2} x, \qquad x \in \text{im}(Q_1^{1/2}),$$

and an unbounded operator

$$Bx(t) = B_t x, \qquad t \in [0, 1),$$

taking values in $C([0, 1), H)$.

The following lemma is crucial for the rest of the paper. Let us recall that the operator $V : H \to E$, where $E$ is a Banach space, is said to be $\gamma$-radonifying if it transforms any cylindrical Gaussian measure on $H$ into a Radon Gaussian measure on $E$.

LEMMA 3.2. (a) *For every* $t \in [0, 1)$, *the operator* $B_t$ *with the domain* $\text{dom}(B) = Q_1^{1/2}(H)$ *extends to a Hilbert–Schmidt operator* $B_t : H \to H$ *and*

$$(3.7) \qquad \int_0^1 \|B_t\|_{\text{HS}} \, dt < \infty.$$

(b) *The operator* $B$ *with the domain* $\text{dom}(B) = Q_1^{1/2}(H)$ *extends to a bounded operator* $B : H \to L^2(0, 1; H)$ *and*

$$(3.8) \qquad |Bx|_{L^2(0,1;H)} = |x|_H, \qquad x \in H.$$

(c) *We have* $\mathcal{K} = \mathcal{L}B$ *and the operator* $\mathcal{K} : H \to C(0, 1; E)$ *is $\gamma$-radonifying.*

PROOF. (a) Note first that $\|Q_1^{-1/2} Q_{1-t}^{1/2}\| \leq 1$ and by (2.4), the operator $Q_{1-t}^{-1/2} S_{1-t}$ is bounded. Therefore, the operator

$$(Q_1^{-1/2} Q_{1-t}^{1/2} Q_{1-t}^{-1/2} S_{1-t} Q^{1/2})^* = \overline{Q^{1/2} S^*_{1-t} Q_1^{-1/2}}$$

is bounded. Moreover, taking (2.7) and (2.6) into account, we obtain, for a certain $C > 0$,

$$\int_0^1 \|B_t^*\|_{\text{HS}} \, dt \leq C \int_0^1 \|Q_t^{-1/2} S_t Q^{1/2}\|_{\text{HS}} \, dt < \infty$$

and (3.7) follows.

For any $h \in H$, we have

$$|Q_1^{1/2} h|^2 = \int_0^1 |Q^{1/2} S^*_{1-t} h|^2 \, dt,$$

and therefore, for $h = Q_1^{-1/2} x$, we obtain

$$(3.9) \qquad |x|^2 = \int_0^1 |Bx(t)|^2 \, dt.$$

Using the density of $Q_1^{1/2}(H)$ in $H$, we can extend (3.9) to the whole of $H$ and (3.8) follows.

(c) For $x \in Q_1^{1/2}(H)$, we have

$$
\begin{aligned}
K_t x &= Q_t S_{1-t}^* Q_1^{-1/2} x \\
&= \int_0^t S_{t-s} Q S_{t-s}^* S_{1-t}^* Q_1^{-1/2} x \, ds \\
&= \int_0^t S_{t-s} Q S_{1-s}^* Q_1^{-1/2} x \, ds \\
&= \int_0^t S_{t-s} Q^{1/2} Bx(s) \, ds = \mathcal{L}(Bx)(t),
\end{aligned}
$$
(3.10)

for all $t \in [0,1]$. By (b), this identity extends to all $x \in H$ and we find that $\mathcal{K} = \mathcal{L}B$ on $H$. By Hypothesis 2.2, we have $\mu(C(0,1;E)) = 1$, hence, $\mathcal{L}: L^2(0,1;H) \to C(0,1;E)$ is $\gamma$-radonifying and therefore, $\mathcal{K} = \mathcal{L}B : H \to C(0,1;E)$ is $\gamma$-radonifying as well. $\square$

We have left open the question of effective verification of Hypothesis 2.3. This is addressed in the following lemma.

LEMMA 3.3. *Assume that either:*

(i) $\dim(H) < \infty$ *or*
(ii) *there exist* $\alpha \in (0,1)$ *and* $\beta < \frac{1+\alpha}{2}$ *such that*

$$\int_0^1 t^{-\alpha} \|S_t Q^{1/2}\|_{\mathrm{HS}}^2 \, dt < \infty$$
(3.11)

*and*

$$\|Q_t^{-1/2} S_t\| \leq \frac{c}{t^\beta}.$$
(3.12)

*Then Hypothesis 2.3 is satisfied.*

PROOF. The proof of (i) extends a classical controllability result from [38] and may be found in [25].

Assume that (ii) holds. By Hypothesis 2.1,

$$Q_t^{-1/2} S_t Q^{1/2} = (Q_t^{-1/2} Q_{t/2}^{1/2})(Q_{t/2}^{-1/2} S_{t/2}) S_{t/2} Q^{1/2},$$

where $\|Q_t^{-1/2} Q_{t/2}^{1/2}\| \leq 1$ and thereby,

$$\|Q_t^{-1/2} S_t Q^{1/2}\|_{\mathrm{HS}} \leq \|Q_{t/2}^{-1/2} S_{t/2}\| \|S_{t/2} Q^{1/2}\|_{\mathrm{HS}}.$$

Therefore, for a certain $c_1 > 0$,

$$\int_0^1 \|Q_t^{-1/2} S_t Q^{1/2}\|_{\text{HS}} \, dt$$
$$\leq \left(\int_0^1 t^\alpha \|Q_{t/2}^{-1/2} S_{t/2}\|^2 \, dt\right)^{1/2} \left(\int_0^1 t^{-\alpha} \|S_{t/2} Q^{1/2}\|_{\text{HS}}^2 \, dt\right)^{1/2}$$
$$\leq c_1 \left(\int_0^1 \frac{1}{t^{2\beta-\alpha}} \, dt\right)^{1/2} \left(\int_0^1 t^{-\alpha} \|S_t Q^{1/2}\|_{\text{HS}}^2 \, dt\right)^{1/2},$$

and (2.7) follows. □

REMARK 3.4. Conditions (3.11) and (3.12) are well known and often used in the theory of SPDE's. Condition (3.11) is a standard assumption that implies the existence of an $H$-continuous version of the OU process $(Z_t^x)$, while (3.12) is closely related to the existence and integrability of the gradient of the OU transition semigroup (cf. [9] for details). Hypothesis 2.3 will be checked in more specific cases in Section 9.

**4. Ornstein–Uhlenbeck bridge.** In Lemma A.2 [applied with $H_1 = L^2(0, 1; H), T = \mathcal{K}$ and $C = Q_1$] an extension of the operator $K_t Q_1^{-1/2}$ to a measurable set $\mathcal{M} \subset H, \mu_1(\mathcal{M}) = 1$ is defined. We use the notation $K_t Q_1^{-1/2}$ for this extension in the present section. Note first that $Q_1^{-1/2} Z_1$ is a cylindrical Gaussian random variable on $H$ and therefore, by Lemma 3.2(b), the process $(K_t Q_1^{-1/2} Z_1)$ is well defined and has $E$-continuous modification. Therefore, we can define an $E$-valued process

$$\widehat{Z}_t = Z_t - K_t Q_1^{-1/2} Z_1, \qquad t \in [0, 1), \quad \text{and} \quad \widehat{Z}_1 = 0,$$

which has an $E$-continuous modification for $t < 1$.

PROPOSITION 4.1. (a) *The $H$-valued Gaussian process $(\widehat{Z}_t)$ is independent of $Z_1$.*

(b) *The covariance operator $\widehat{Q}_t$ of $\widehat{Z}_t$ is given by*

(4.1) $$\widehat{Q}_t = Q_t^{1/2}(I - V_t^* V_t) Q_t^{1/2}, \qquad t \in [0, 1).$$

(c) *The process $(\widehat{Z}_t)$ is continuous in $E$ for $t \in [0, 1]$.*

PROOF. (a) For $h, k \in \operatorname{im}(Q_1^{1/2})$, (3.1) yields

$$\mathbb{E}\langle \widehat{Z}_t, h\rangle\langle Z_1, k\rangle = \mathbb{E}\langle Z_t, h\rangle\langle Z_1, k\rangle - \mathbb{E}\langle K_t Q_1^{-1/2} Z_1, h\rangle\langle Z_1, k\rangle$$
$$= \langle S_{1-t} Q_t h, k\rangle - \langle Q_1 Q_1^{-1/2} K_t^* h, k\rangle$$
$$= \langle S_{1-t} Q_t h, k\rangle - \langle S_{1-t} Q_t h, k\rangle = 0$$

and therefore, the process $(\widehat{Z}_t)$ and $Z_1$ are independent.

(b) It follows from (a) that
$$Q_t = \widehat{Q}_t + K_t K_t^*.$$
Hence, the definition of $K_t$ and $V_t$ yields
$$\widehat{Q}_t = Q_t - Q_t S_{1-t}^* Q_1^{-1} S_{1-t} Q_t = Q_t^{1/2}(I - V_t^* V_t) Q_t^{1/2}.$$

(c) Using (3.2), we find easily that
$$\lim_{t \to 0} \operatorname{tr}(\widehat{Q}_t) = 0. \tag{4.2}$$

To prove that
$$\lim_{t \to 1} \operatorname{tr}(\widehat{Q}_t) = 0, \tag{4.3}$$

we note first that
$$\operatorname{tr}(\widehat{Q}_t) = \operatorname{tr}((I - V_t^* V_t)(Q_t - Q_1)) + \operatorname{tr}((I - V_t^* V_t) Q_1).$$

Next, it is easy to see that
$$0 \le \lim_{t \to 1} \operatorname{tr}((I - V_t^* V_t)(Q_1 - Q_t)) \le \lim_{t \to 1} \operatorname{tr}(Q_1 - Q_t) = 0.$$

We have also
$$\operatorname{tr}((I - V_t^* V_t) Q_1) = \operatorname{tr}(Q_1) - \sum_{k=1}^{\infty} |V_t Q_1^{1/2} e_k|^2,$$

and (4.3) follows from (3.3), (3.2) and the dominated convergence. Since the process $(\widehat{Z}_t)$ has $E$-continuous version by Lemma 3.2(b), it follows that
$$0 = \lim_{t \to 1} \widehat{Z}_t = \widehat{Z}_1.$$

This fact completes the proof of continuity. $\square$

PROPOSITION 4.2. *There exists a Borel subspace* $\mathcal{M} \subset H$ *such that* $\mu_1(\mathcal{M}) = 1$ *and for all* $x \in H$ *and* $y \in \mathcal{M}$, *the $H$-valued Gaussian process*
$$\widehat{Z}_t^{x,y} = Z_t^x - K_t Q_1^{-1/2}(Z_1^x - y) \tag{4.4}$$
*is well defined for all $t \in [0, 1)$. Moreover,*
$$\widehat{Z}_t^{x,y} = S_t x - K_t Q_1^{-1/2} S_1 x + K_t Q_1^{-1/2} y + \widehat{Z}_t, \qquad \mathbb{P}\text{-a.s.} \tag{4.5}$$

PROOF. By Lemma A.3, we can choose a measurable linear space $\mathcal{M}$ such that $K_t Q_1^{-1/2}$ is linear on $\mathcal{M}$ with $\mu_1(\mathcal{M}) = 1$ and the mapping $(t,y) \to K_t Q_1^{-1/2} y$ is measurable. By (2.4) we have $S_1 x \in \mathrm{im}(Q_1^{1/2})$ and therefore, $K_t Q_1^{-1/2}(Z_1^x - y)$ is well defined for any $y \in \mathcal{M}$. Clearly, $\widehat{Z}_t^{x,y}$ may be rewritten in the form (4.5). $\square$

The process $(\widehat{Z}_t^{x,y})$ defined in Proposition 4.2 will be called an Ornstein–Uhlenbeck bridge on $H$ (connecting points $x \in H$ and $y \in \mathcal{M}$). We will denote by $\widehat{\mu}^{x,y}$ the law of the process $\{\widehat{Z}_t^{x,y} : t \in [0,1]\}$.

THEOREM 4.3. *There exists a Borel subspace $\mathcal{M} \subset E$ with $\mu_1(\mathcal{M}) = 1$, such that the process $(\widehat{Z}_t^{x,y})$ has $E$-continuous version for each $x \in E$ and $y \in \mathcal{M}$. Moreover, there exists a measurable mapping $U : \mathcal{M} \to \mathbb{R}_+$ and a random variable $k$, such that*

$$(4.6) \qquad \|\widehat{Z}^{x,y}\|_{C(0,1;E)} \leq k(1 + |x|_E + U(y)), \qquad x \in E, y \in \mathcal{M},$$

*and*

$$(4.7) \qquad \mathbb{E}\|\widehat{Z}^{x,y}\|_{C(0,1;E)}^n \leq L(n)(1 + |x|_E^n + (U(y))^n)$$

*for each $n \in \mathbb{N}$, $x \in E$ and $y \in \mathcal{M}$, where $L(n)$ is a constant depending on $n$ only.*

PROOF. It was already shown in Proposition 4.1 that the process $(\widehat{Z}_t)$ has trajectories in $C(0,1;E)$ and we have

$$(4.8) \qquad k_1 = \sup_{t \leq 1} |\widehat{Z}_t|_E < \infty, \qquad \mathbb{P}\text{-a.s.}$$

Since $(\widehat{Z}_t)$ is a Gaussian process, we obtain, for any $m > 0$,

$$(4.9) \qquad k_2(m) = \mathbb{E} \sup_{t \leq 1} |\widehat{Z}_t|_E^m < \infty.$$

The same argument shows that the process $t \to K_t Q_1^{-1/2} y$ has trajectories in $C(0,1;E)$ for every $y \in \mathcal{M}$, where $\mathcal{M}$ is given by Proposition 4.2 (possibly, after excluding a zero $\mu_1$-measure set). By the strong Feller property we have $S_1 x \in \mathrm{im}(Q_1^{1/2}) \subset \mathcal{M}$, so it follows from Proposition 4.2 that $\widehat{\mu}^{x,y}(C(0,1;E)) = 1$ for $x \in E$ and $y \in \mathcal{M}$. Furthermore, using the notation from Lemma 3.2, we have $K_t Q_1^{-1/2} S_1 x = \mathcal{K}(Q_1^{-1/2} S_1 x)(t), t \in [0,1]$. By (2.4), the operator $Q_1^{-1/2} S_1$ is bounded and it is easy to see that (3.3) together with (3.6) yields continuity of the mapping $K_t Q_1^{-1/2} S_1 : E \to C(0,1;E)$. Hence, setting $U(y) = \|\mathcal{K} Q_1^{-1/2} y\|_{C(0,1;E)}$ and taking into account (4.8) and (4.9), we obtain both (4.6) and (4.7) for all $n > 0$. $\square$

The following theorem justifies the intuitive notion of the OU bridge $(\widehat{Z}_t^{x,y})$ given in the Introduction.

THEOREM 4.4. *Let $\Phi:C(0,1;E) \to \mathbb{R}$ be a Borel mapping such that, for $x \in E$,*

$$\mathbb{E}|\Phi(Z^x)| < \infty.$$

*Then*

$$\mathbb{E}(\Phi(Z^x)|Z_1^x = y) = \mathbb{E}\Phi(\widehat{Z}^{x,y}), \qquad \mu_1\text{-}a.e.$$

PROOF. By Hypothesis 2.2 and Theorem 4.3, the processes $(Z_t^x)$ and $(\widehat{Z}_t^{x,y})$ are concentrated on $C(0,1;E)$ and $(K_t Q_1^{-1/2}(S_1 x - y)) \in E$. Moreover, the processes $(\widehat{Z}_t^{x,y})$ and $(K_t Q_1^{-1/2} Z_1^x)$ are independent by Proposition 4.1. Therefore, using well-known properties of conditional expectations, we obtain

$$\mathbb{E}(\Phi(Z^x)|Z_1^x = y) = \mathbb{E}(\Phi(\widehat{Z}^{x,y} + K_t Q_1^{-1/2}(Z_1^x - y))|Z_1^x = y)$$
$$= \mathbb{E}\Phi(\widehat{Z}^{x,y}), \qquad \mu_1\text{-a.e.}$$

for any $x \in E$. □

Let

$$Y_u = \int_u^1 S_{1-s} Q^{1/2}\, dW_s, \qquad u \leq 1,$$

and

$$H_u = Q_{1-u}^{-1/2} S_{1-u} Q^{1/2}, \qquad u < 1.$$

LEMMA 4.5. *For all $u \in [0,1]$, we have*

(4.10) $$Y_u = Q_{1-s} Q_1^{-1} Z_1 - S_{1-s} \widehat{Z}_s, \qquad \mathbb{P}\text{-}a.s.,$$

*where $Q_{1-s} Q_1^{-1}$ is bounded for all $s \in [0,1]$.*

PROOF. We have

$$K_t Q_1^{-1/2} Z_1 = \left(\int_0^t S_{t-s} Q^{1/2} H_s^* \, ds\right) Q_1^{-1/2} Z_1$$

and

$$S_{1-t} K_t Q_1^{-1/2} Z_1 = \left(\int_0^t S_{1-s} Q^{1/2} H_s^* \, ds\right) Q_1^{-1/2} Z_1$$
$$= (Q_1 - Q_{1-t}) Q_1^{-1} Z_1$$
$$= Z_1 - Q_{1-t} Q_1^{-1} Z_1,$$

and thereby,
$$Z_1 - S_{1-t}K_t Q_1^{-1/2} Z_1 = Q_{1-t}Q_1^{-1}Z_1.$$

Therefore, by definition of $\widehat{Z}_t$, we obtain
$$\begin{aligned} Y_s &= Z_1 - S_{1-s}Z_s \\ &= Z_1 - S_{1-s}(Z_s - K_s Q_1^{-1/2}Z_1) - S_{1-s}K_s Q_1^{-1/2}Z_1 \\ &= Z_1 - S_{1-s}\widehat{Z}_s - (Z_1 - Q_{1-s}Q_1^{-1}Z_1) \\ &= Q_{1-s}Q_1^{-1}Z_1 - S_{1-s}\widehat{Z}_s. \quad \square \end{aligned}$$

Since the operator-valued function $t \to Q_t$ is continuous in the weak operator topology and all the operators $Q_t$ are compact for $t > 0$, there exists a measurable choice of eigenvectors $\{e_k(t): k \geq 1\}$ and eigenvalues $\{\lambda_k(t): k \geq 1\}$. For each $n \geq 1$ we define a process
$$\alpha_u^n = \sum_{k=1}^n \frac{1}{\sqrt{\lambda_k(1-u)}} \langle Y_u, e_k(1-u)\rangle H_u^* e_k(1-u).$$

LEMMA 4.6.   *There exists a measurable stochastic process $(\alpha_u)$ defined on $[0,1)$ such that, for each $a < 1$,*

(4.11) $$\lim_{n\to\infty} \mathbb{E} \int_0^a |\alpha_u^n - \alpha_u|^2 \, du = 0$$

*and for each $h \in H$ and $a < 1$, the series*

(4.12) $$\langle \alpha_u, h \rangle = \sum_{k=1}^\infty \frac{1}{\sqrt{\lambda_k(1-u)}} \langle Y_u, e_k(1-u)\rangle \langle e_k(1-u), H_u h\rangle$$

*converges in $L^2(0,a)$ in mean square. Moreover, if $0 \leq u \leq v < 1$, then, for all $h, k \in H$,*

(4.13) $$\mathbb{E}\langle \alpha_u, h\rangle \langle \alpha_v, k\rangle = \langle H_u h, Q_{1-u}^{-1/2} Q_{1-v}^{1/2} H_v k\rangle,$$

*where the operator $Q_{1-u}^{-1/2}Q_{1-v}^{1/2}$ is bounded. Finally,*

(4.14) $$\mathbb{E} \int_0^1 |\alpha_u| \, du < \infty.$$

In what follows we will use the notation $\langle \alpha_u, h\rangle = \langle Q_{1-u}^{-1/2} Y_u, H_u h\rangle$.

PROOF OF LEMMA 4.6.   For $u \leq v \leq 1$,

(4.15) $$\mathbb{E}\langle Y_u, h\rangle \langle Y_v, k\rangle = \langle Q_{1-v}h, k\rangle, \qquad h, k \in H.$$

Therefore,

$$\mathbb{E}\langle \alpha_u^n - \alpha_u^m, h\rangle^2 = \sum_{j=m+1}^{n} \frac{1}{\lambda_k(1-u)} \mathbb{E}\langle Y_u, e_k(1-u)\rangle^2 \langle e_k(1-u), H_u h\rangle^2$$

(4.16)

$$= \sum_{j=m+1}^{n} \langle e_k(1-u), H_u h\rangle^2 \underset{n,m\to\infty}{\longrightarrow} 0,$$

hence, the process

$$\langle \alpha_u, h\rangle = \sum_{k=1}^{\infty} \frac{1}{\sqrt{\lambda_k(1-u)}} \langle Y_u, e_k(1-u)\rangle \langle e_k(1-u), H_u h\rangle$$

(4.17)

$$= \langle Q_{1-u}^{-1/2} Y_u, H_u h\rangle$$

is well defined and measurable for each $h \in H$ and $u < 1$. Let $P_n$ be an orthogonal projection on $\lin\{e_k(1-v) : k \leq n\}$ and $H_u^n = P_n H_u$. Then $Q_{1-u}^{-1/2} H_u^n$ is bounded on $H$ and we may define $\alpha_u^n = (Q_{1-u}^{-1/2} H_u^n)^* Y_u$. By (4.15),

$$\mathbb{E}\langle \alpha_u^n, h\rangle \langle \alpha_v^n k\rangle = \langle Q_{1-v} Q_{1-u}^{-1/2} H_u^n h, Q_{1-v}^{-1/2} H_v^n k\rangle$$

$$= \langle H_u^n h, Q_{1-u}^{-1/2} Q_{1-v}^{1/2} H_v^n k\rangle.$$

By (2.6), the operator $Q_{1-u}^{-1/2} Q_{1-v}^{1/2}$ is bounded and, therefore,

$$\mathbb{E}\langle Q_{1-u}^{-1/2} Y_u, H_u h\rangle \langle Q_{1-v}^{-1/2} Y_v, H_v k\rangle = \lim_{n\to\infty} \mathbb{E}\langle \alpha_u^n, h\rangle \langle \alpha_v^n, k\rangle$$

$$= \langle H_u h, Q_{1-u}^{-1/2} Q_{1-v}^{1/2} H_v k\rangle.$$

It follows from (4.13) that

$$\mathbb{E}\langle \alpha_u, h\rangle^2 = |H_u h|^2,$$

hence,

$$\mathbb{E}|\alpha_u|^2 = \|H_u\|_{\mathrm{HS}}^2 < \infty, \qquad u < 1,$$

and by Hypothesis 2.3,

$$\mathbb{E}\int_0^1 |\alpha_u|\, du < \infty.$$

Then (4.16) and the dominated convergence yield

$$\lim_{n,m\to\infty} \int_0^a \mathbb{E}|\alpha_u^n - \alpha_u^m|^2\, du = 0.$$

As a consequence, we find that (4.11) holds for any $a \in (0,1)$. $\square$

LEMMA 4.7. *The cylindrical process*

$$\zeta_t = W_t - \int_0^t \alpha_u \, du, \qquad t \leq 1,$$

*is a standard cylindrical Wiener process on $H$ independent of $Z_1$.*

PROOF. We need to show that, for any $h \in H$, the process

$$\langle \zeta_t, h \rangle = \langle W_t, h \rangle - \int_0^t \langle Q_{1-u}^{-1/2} Y_u, H_u h \rangle$$

is a real-valued Wiener process. Let $h, k \in H$. We will show first that, for $r < t < 1$,

(4.18) $$\mathbb{E}\langle \zeta_t - \zeta_r, h \rangle \langle \zeta_r, k \rangle = 0.$$

We have

$$\mathbb{E}\langle \zeta_t - \zeta_r, h \rangle \langle \zeta_r, k \rangle$$
$$= -\mathbb{E}\langle W_t - W_r, h \rangle \int_0^r \langle Q_{1-u}^{-1/2} Y_u, H_u k \rangle \, du$$
$$\quad - \mathbb{E}\langle W_r, k \rangle \int_r^t \langle Q_{1-u}^{-1/2} Y_u, H_u h \rangle \, du$$
$$\quad + \mathbb{E}\left(\int_r^t \langle Q_{1-u}^{-1/2} Y_u, H_u h \rangle \, du\right)\left(\int_0^r \langle Q_{1-u}^{-1/2} Y_u, H_v k \rangle \, dv\right)$$
$$= -I_1 - I_2 + I_3.$$

We will consider $I_1$ first. Taking into account that the series (4.17) is mean-square convergent, for each $u \in (0, 1)$, we obtain

$$\mathbb{E}\langle W_t - W_r, h \rangle \langle Q_{1-u}^{-1/2} Y_u, H_u k \rangle$$
$$= \sum_{n=1}^{\infty} \frac{\langle e_n(1-u), H_u k \rangle}{\sqrt{\lambda_n(1-u)}} \mathbb{E}(\langle Y_u, e_n(1-u) \rangle \langle W_t - W_r, h \rangle).$$

Next, for $u \leq r$,

$$\mathbb{E}(\langle Y_u, e_n(1-u) \rangle \langle W_t - W_r, h \rangle)$$
$$= \mathbb{E} \int_u^1 \langle Q^{1/2} S_{1-s}^* e_n(1-u), dW_s \rangle \int_r^t \langle h, dW_s \rangle$$
$$= \int_r^t \langle S_{1-s} Q^{1/2} h, e_n(1-u) \rangle \, ds,$$

and therefore,

$$\mathbb{E}\langle W_t - W_r, h\rangle\langle Q_{1-u}^{-1/2}Y_u, H_u k\rangle$$
$$= \sum_{n=1}^{\infty} \frac{\langle e_n(1-u), H_u k\rangle}{\sqrt{\lambda_n(1-u)}} \int_r^t \langle S_{1-s}Q^{1/2}h, e_n(1-u)\rangle \, ds$$
$$= \int_r^t \left(\sum_{n=1}^{\infty} \frac{\langle e_n(1-u), H_u k\rangle}{\sqrt{\lambda_n(1-u)}} \langle S_{1-s}Q^{1/2}h, e_n(1-u)\rangle \, ds\right)$$
$$= \int_r^t \langle Q_{1-u}^{-1/2} S_{1-s} Q^{1/2} h, H_u k\rangle \, ds$$
$$= \int_r^t \langle Q_{1-u}^{-1/2} Q_{1-s}^{1/2} H_s h, H_u k\rangle \, ds$$

and

$$(4.19) \qquad I_1 = \int_0^r \int_r^t \langle Q_{1-u}^{-1/2} Q_{1-s}^{1/2} H_s h, H_u k\rangle \, ds \, du,$$

where the operator $Q_{1-u}^{-1/2} Q_{1-s}^{1/2}$ is bounded. By similar arguments, we find that $I_2 = 0$ and (4.13) yields

$$I_3 = \int_0^r \int_r^t \langle Q_{1-v}^{-1/2} Q_{1-u}^{1/2} H_u h, H_v k\rangle \, du \, dv$$

and in view of (4.19), $I_1 = I_3$, and since $I_2 = 0$, (4.18) follows. We will show, that for $h \in H$ and $t < 1$,

$$(4.20) \qquad \mathbb{E}\langle \zeta_t, h\rangle^2 = t|h|^2.$$

We have

$$\mathbb{E}\langle \zeta_t, h\rangle^2 = t|h|^2 - 2\mathbb{E}\langle W_t, h\rangle \int_0^t \langle Q_{1-u}^{-1/2} Y_u, H_u h\rangle$$
$$+ \mathbb{E}\left(\int_0^t \langle Q_{1-u}^{-1/2} Y_u, H_u h\rangle \, du\right)\left(\int_0^t \langle Q_{1-v}^{-1/2} H_v Y_v, k\rangle \, dv\right)$$
$$= t|h|^2 - 2J_1 + J_3.$$

Proceeding in the same way as in the computation of $I_1$ and $I_3$, we obtain

$$(4.21) \qquad J_1 = \int_0^t \mathbb{E}\langle W_t, h\rangle\langle Q_{1-u}^{-1/2} Y_u, H_u h\rangle \, du$$
$$= \int_0^t \int_u^t \langle Q_{1-u}^{-1/2} Q_{1-s}^{1/2} H_s h, H_u h\rangle \, ds \, du.$$

Invoking again (4.13), we obtain

$$J_3 = \int_0^t \int_0^t \mathbb{E}\langle Q_{1-u}^{-1/2} Y_u, H_u h\rangle \langle Q_{1-v}^{-1/2} Y_v, H_v h\rangle \, du \, dv$$

$$= \int_0^t \int_0^v \langle H_u h, Q_{1-u}^{-1/2} Q_{1-v}^{1/2} H_v h\rangle \, du \, dv$$

$$+ \int_0^t \int_v^t \langle Q_{1-v}^{-1/2} Q_{1-u}^{1/2} H_u h, H_v h\rangle \, du \, dv.$$

Since $t < 1$ and the functions under the integrals are continuous, we can change the order of integration in the first integral and obtain

$$J_3 = \int_0^t \int_u^t \langle H_u h, Q_{1-u}^{-1/2} Q_{1-v}^{1/2} H_v h\rangle \, dv \, du$$

$$+ \int_0^t \int_v^t \langle Q_{1-v}^{-1/2} Q_{1-u}^{1/2} H_u h, H_v h\rangle \, du \, dv.$$

Hence, $J_3 = 2J_1$ and (4.20) follows. Combining (4.18) and (4.20), we find that, for $s, t < 1$,

$$\mathbb{E}\langle \zeta_s, h\rangle \langle \zeta_t, k\rangle = \min(s,t)\langle h, k\rangle, \qquad h, k \in H.$$

Since

$$\sup_{t<1} \mathbb{E}|\langle \zeta_t, h\rangle| = |h|^2, \qquad h \in H,$$

there exists a cylindrical random variable $\zeta_1$ such that

$$\lim_{t \to 1} \langle \zeta_t, h\rangle = \langle \zeta_1, h\rangle, \qquad h \in H,$$

for all $t \leq 1$. Therefore, $(\zeta_t)$ is a cylindrical Brownian motion for $t \in [0,1]$. It remains to show that, for any $t < 1$,

(4.22) $$\mathbb{E}\langle \zeta_t, h\rangle \langle Z_1, k\rangle = 0, \qquad h, k \in H.$$

By definition of $\zeta_t$, it is enough to show that

(4.23) $$\mathbb{E}\langle W_t, h\rangle \langle Z_1, k\rangle = \mathbb{E} \int_0^t \langle \alpha_u, h\rangle \langle Z_1, k\rangle \, du.$$

Now, we have

(4.24) $$\mathbb{E}\langle W_t, h\rangle \langle Z_1, k\rangle = \int_0^t \langle S_{1-u} Q^{1/2} h, k\rangle \, du.$$

Invoking (4.10) and using the fact that $Z_1$ and $(\widehat{Z}_t)$ are independent, we obtain

$$\mathbb{E} \int_0^t \langle \alpha_u, h \rangle \langle Z_1, k \rangle \, du$$

$$= \int_0^t \mathbb{E} \langle Q_{1-u}^{-1/2} Y_u, H_u h \rangle \langle Z_1, k \rangle \, du,$$

(4.25)
$$\int_0^t \mathbb{E} \langle Q_{1-u}^{-1/2} Q_{1-u} Q_1^{-1} Z_1, H_u h \rangle \langle Z_1, k \rangle \, du$$

$$= \int_0^t \langle Q_{1-u}^{1/2} H_u h, k \rangle \, du$$

$$= \int_0^t \langle S_{1-u} Q^{1/2} h, k \rangle \, du.$$

Comparing (4.24) and (4.25), we obtain (4.23) and the lemma follows. □

REMARK 4.8. Lemma 4.7 allows to define the Ornstein–Uhlenbeck bridge $(\widehat{Z}_t^{x,y})$ as a unique solution of a certain linear stochastic evolution equation. As it is not needed in this paper, it is omitted, see [15] for details.

PROPOSITION 4.9. *Let*

$$B_1(s) = (Q_{1-s}^{-1/2} S_{1-s} Q^{1/2})^* Q_{1-s}^{-1/2} S_{1-s},$$

$$B_2(s) = (Q_1^{-1/2} S_{1-s} Q^{1/2})^* Q_1^{-1/2} S_1,$$

$$B_3(s) y = (Q_1^{-1/2} S_{1-s} Q^{1/2})^* Q_1^{-1/2} y, \qquad y \in \mathrm{im}(Q_1^{1/2}), s \in (0,1).$$

*Then*

(4.26) $$\mathbb{E} \int_0^1 |B_1(s) \widehat{Z}_s| \, ds < \infty,$$

(4.27) $$\int_0^1 |B_2(s) x|^2 \, ds = |Q_1^{-1/2} S_1 x|^2, \qquad x \in H.$$

*Moreover, there exists a Borel subspace* $\mathcal{M} \subset H$ *with* $\mu_1(\mathcal{M}) = 1$ *such that* $B_3$ *extends to a linear mapping* $B_3 \colon \mathcal{M} \to L^1(0,1;H)$, *that is,*

(4.28) $$\int_0^1 |B_3(s) y| \, ds < \infty, \qquad y \in \mathcal{M},$$

*and* $B_3(s) Z_1$ *has the covariance* $(Q_1^{-1/2} S_{1-s} Q^{1/2})^* (Q_1^{-1/2} S_{1-s} Q^{1/2})$ *for each* $s \in [0,1)$.

PROOF. Condition (4.26) is a reformulation of (4.14) in Lemma 4.6, where it was also shown that the function $t \to B_1(t)\widehat{Z}_t$ is integrable. Invoking the definition of the operator $B$ and Lemma 3.2, we obtain

$$(4.29) \quad \int_0^1 |B_2(s)x|^2\,ds = \int_0^1 |BQ_1^{-1/2}S_1 x(s)|^2\,ds = |Q_1^{-1/2}S_1 x|^2,$$

and (4.27) follows.

For $y \in \mathrm{im}(Q_1^{1/2})$, we have $B_3(t)y = BQ_1^{-1/2}y$ and for any $a < 1$,

$$\int_0^a \|B_t\|_{\mathrm{HS}}^2\,dt < \infty,$$

and taking (3.7) into account, we may apply Lemma A.3, which yields the desired result. $\square$

**5. Transition density of semilinear stochastic evolution equation.** In this section a formula for transition densities defined by equation (2.1) will be derived and some useful lower estimates on transition densities will be established.

PROPOSITION 5.1. *Assume that equation* (2.1) *has an invariant measure* $\mu^* \in \mathcal{P}$. *Then*

$$(5.1) \quad \|P_t^*\nu - \mu^*\|_{\mathrm{var}} \to 0, \qquad t \to \infty, \nu \in \mathbb{P},$$

*where* $\|\cdot\|_{\mathrm{var}}$ *denotes the total variation of measures. Furthermore, for each* $x \in E$ *and* $T > 0$, *the measures* $P(T, x, \cdot)$ *and* $\mu_T^x$ *are equivalent and for* $\mu_T$ *a.e.* $y$,

$$(5.2) \quad \frac{dP(T, x, \cdot)}{d\mu_T^x}(y) = \mathbb{E}\left(\exp\left(\int_0^T \langle G(Z_t^x), dW_t\rangle - \frac{1}{2}\int_0^T |G(Z_t^x)|^2\,dt\right)\Big| Z_T^x = y\right).$$

PROOF. By Hypothesis 2.1, the Ornstein–Uhlenbeck process is strongly Feller, therefore, its distributions $(\mu_t^x)$ are equivalent for $x \in H$, $t > 0$. If (5.2) holds for each $T > 0$ and $x \in E$, we have the equivalence $P(T, x, \cdot) \sim \mu_T^x$, hence, $(P(T, x, \cdot))_{T>0, x\in E}$ are equivalent, as well and the convergence (5.1) follows from well-known results (cf. [37, 40]).

It remains to prove (5.2), which follows from the Girsanov theorem (see, e.g., [9], Theorem 10.4). In order to apply this result, we must verify (taking, for simplicity, $T = 1$)

$$(5.3) \quad \mathbb{E}\exp\rho(Z^x) = 1,$$

where

$$\rho(Z^x) := \int_0^1 \langle G(Z_s^x), dW_s \rangle - \tfrac{1}{2} \int_0^1 |G(Z_s^x)|^2 \, ds. \tag{5.4}$$

For $n \geq 1$, set

$$F_n(x) = \begin{cases} F(x), & \text{if } |x|_E \leq n, \\ F\left(\dfrac{nx}{|x|_E}\right), & \text{if } |x|_E > n, \end{cases} \tag{5.5}$$

and let $G_n$ be defined by $F_n(x) := Q^{1/2} G_n(x)$. Obviously the approximating equations

$$\begin{aligned} dX_n(t) &= (AX_n(t) + F_n(X_n(t))) + \sqrt{Q}\, dW_t, \\ X_n(0) &= x, \end{aligned} \tag{5.6}$$

have uniquely defined solutions $\mathbb{P}$-a.s. in $C(0,1;E)$ and denoting by $\tilde{P}_X, \tilde{P}_{X_n}$ and $\tilde{P}_{Z^x}$ the distributions in $C(0,1;E)$ of $X$, $X_n$ and $Z^x$, respectively, we have

$$\lim_{n \to \infty} \mathbb{P}\left( \sup_{t \in [0,1]} |X_n(t) - X(t)|_E > 0 \right) = 0. \tag{5.7}$$

Hence,

$$\|\tilde{P}_{X_n}(\cdot) - \tilde{P}_X\|_{\text{var}} \to 0, \qquad m \to \infty, \tag{5.8}$$

thus, $(\tilde{P}_{X_n})$ is a Cauchy sequence in the metric of total variation. Therefore, the sequence of densities

$$\frac{d\tilde{P}_{X_n}}{d\tilde{P}_{Z^x}} = \exp \rho_n(Z^x), \tag{5.9}$$

where

$$\rho_n(Z^x) := \int_0^1 \langle G_n(Z_s^x), dW_s \rangle - \tfrac{1}{2} \int_0^1 |G_n(Z_s^x)|^2 \, ds, \tag{5.10}$$

is a Cauchy, hence, convergent, sequence in $L^1(\Omega)$. As $G_n$ is bounded for each $n$, obviously $\mathbb{E} \exp \rho_n(Z^x) = 1$, so it remains to identify the $L^1(\Omega)$-limit of $\exp \rho_n$ with $\exp \rho$. Clearly, $G_n \to G$ pointwise and $|G_n(x)| \leq K(1 + |x|_E^m)$, therefore,

$$\mathbb{E} \left| \int_0^1 \langle G_n(Z_s^x), dW_s \rangle - \int_0^1 \langle G(Z_s^x), dW_s \rangle \right|^2$$

$$= \mathbb{E} \int_0^1 |G_n(Z_s^x) - G(Z_s^x)|^2 \, ds \to 0$$

by the dominated convergence theorem and Hypothesis 2.2. Similarly, we have
$$\int_0^1 |G_n(Z_s^x)|^2 \, ds \to \int_0^1 |G(Z_s^x)|^2 \, ds, \qquad \mathbb{P}\text{-a.s.},$$
so we obtain (possibly, for a subsequence) $\exp \rho_n(Z^x) \to \exp \rho(Z^x)$ $\mathbb{P}$-a.s., which completes the proof of (5.2). $\square$

We will now state one of our main results, which provides a formula for the density $\frac{dP_t^*\nu}{d\mu_t}$ for a given time $t > 0$ (we may take $t = 1$). It follows from the Fubini theorem that the density has the form

(5.11)
$$\frac{dP_1^*\nu}{d\mu_1}(y) = \int_E \frac{dP(1,x,\cdot)}{d\mu_1}(y)\nu(dx)$$
$$= \int_E \frac{dP(1,x,\cdot)}{d\mu_1^x}(y)\frac{d\mu_1^x}{d\mu_1}(y)\nu(dx), \qquad \mu_1\text{-a.e.,}$$

provided the product of densities inside the integral on the r.h.s. is $(x,y)$-measurable. As mentioned in the preceding proof, the Gaussian measures $\mu_1^x$ and $\mu_1$ are equivalent with the density given by the Cameron–Martin formula

(5.12)
$$g(x,y) := \frac{d\mu_1^x}{d\mu_1}(y)$$
$$= \exp\Big(\langle Q_1^{-1/2}S_1 x, Q_1^{-1/2} y\rangle - \frac{1}{2}|Q_1^{-1/2}S_1 x|^2\Big), \qquad x \in E,$$

for $\mu_1$-almost all $y \in E$.

THEOREM 5.2. *For each $\nu \in \mathcal{P}$, we have*

(5.13) $$\frac{dP_1^*\nu}{d\mu_1}(y) = \int_E h(x,y)g(x,y)\nu(dx), \qquad \mu_1\text{-a.e.,}$$

*where $g$ is defined by the Cameron–Martin formula* (5.12), *and for $x \in E$ and $\mu_1$-almost all $y \in E$,*

(5.14)
$$h(x,y) := \mathbb{E}\exp\Big(\rho(\widehat{Z}^{x,y})$$
$$\qquad - \int_0^1 \langle G(\widehat{Z}_s^{x,y}), B_1(s)\widehat{Z}_s + B_2(s)x - B_3(s)y\rangle \, ds\Big),$$

*where $B_1$, $B_2$ and $B_3$ are defined in Proposition* 4.9. *In particular, for each $x \in E$, we have*

(5.15) $$\frac{dP(1,x,\cdot)}{d\mu_1}(y) = h(x,y)g(x,y), \qquad \mu_1\text{-a.e.}$$

PROOF. Since both $g$ and $h$ are $(x, y)$-measurable, taking into account (5.11) and (5.12), we only have to prove that

$$\frac{dP(1, x, \cdot)}{d\mu_1^x}(y) = h(x, y), \qquad x \in E, \mu_1\text{-a.e.} \tag{5.16}$$

Assume at first that the mapping $G$ is bounded and let $t_i^k := \frac{i}{k}$ for $k \in \mathbb{N}$, $i = 0, 1, \ldots, k$, that is, $\Delta_k = \{t_0^k, t_1^k, \ldots, t_k^k\}$ are equidistant divisions of the interval $[0, 1]$, $t_0^k = 0$, $t_k^k = 1$, $t_{i+1}^k - t_i^k = 1/k$ (for brevity, the dependence of $t_i^k$ on $k$ will be suppressed in the notation). Set, for $k \geq 1$,

$$\rho^k(Z^x) := \sum_{i=0}^{k-1} \langle G(Z_{t_i}^x), W_{t_{i+1}} - W_{t_i} \rangle - \tfrac{1}{2} \int_0^1 |G(Z_s^x)|^2 \, ds. \tag{5.17}$$

The mapping $G$ is assumed to be bounded, thus,

$$\int_0^1 |G(Z_s^x)|^2 \, ds + \int_0^1 \left| \sum_{i=0}^{k-1} G(Z_{t_i}^x) \mathbf{1}_{[t_i, t_{i+1}]}(s) \right|^2 ds \leq 2 \sup |G|^2 < \infty, \tag{5.18}$$

so the random variables $\exp \rho^k(Z^x)$ are uniformly integrable on $\Omega$. Clearly, $\exp \rho^k(Z^x) \to \exp \rho(Z^x)$ $\mathbb{P}$-a.s. (possibly, for a subsequence), and therefore,

$$\exp \rho(Z^x) = \lim_{k \to \infty} \exp \rho^k(Z^x) \qquad \text{in } L^1(\Omega), \tag{5.19}$$

which in view of (5.2) yields

$$\frac{dP(1, x, \cdot)}{d\mu_1^x}(y) = \mathbb{E}(\exp \rho(Z^x) | Z_1^x = y) = \lim_{k \to \infty} \mathbb{E}(\exp \rho^k(Z^x) | Z_1^x = y) \tag{5.20}$$

for $\mu_1^x$-almost all $y \in H$. On the other hand, in terms of the cylindrical Wiener process $\zeta_t$ defined in Lemma 4.7, we have

$$\rho^k(Z^x) = \sum_{i=0}^{k-1} \langle G(Z_{t_i}^x), \zeta_{t_{i+1}} - \zeta_{t_i} \rangle$$
$$+ \int_{t_i}^{t_{i+1}} \langle G(Z_{t_i}^x), H_s^* Q_{1-s}^{-1/2} Y_s \rangle \, ds - \tfrac{1}{2} \int_0^1 |G(Z_s^x)|^2 \, ds. \tag{5.21}$$

Invoking (4.10), we obtain

$$- H_s^* Q_{1-s}^{-1/2} Y_s = H^* Q_{1-s}^{-1/2} S_{1-s} \widehat{Z}_s - H_s^* Q_{1-s}^{1/2} Q_1^{-1/2} Z_1$$
$$= B_1(s) \widehat{Z}_s - B_3(s) Z_1 \tag{5.22}$$

for $s \in (0, 1)$ $\mathbb{P}$-a.s. [note that both terms on the r.h.s. of (5.22) are well defined $\mathbb{P}$-a.s. in $L^1(0, 1, H)$ by Proposition 4.9]. Therefore, for $y \in \mathcal{M}$, where $\mathcal{M}$ is the intersection of the two full measure sets (denoted in both cases by $\mathcal{M}$) from Theorem 4.3 and Proposition 4.9, respectively, we obtain

$$\mathbb{E}(\exp \rho^k(Z^x)|Z_1^x = y)$$
$$= \mathbb{E}(\exp \rho^k(Z^x)|Z_1 = y - S_1 x)$$
$$= \mathbb{E}\left(\exp\left(\sum_{i=0}^{k-1}\langle G(Z_{t_i}^x), \zeta_{t_{i+1}} - \zeta_{t_i}\rangle\right.\right.$$
$$- \int_{t_i}^{t_{i+1}} \langle G(Z_{t_i}^x), B_1(s)\widehat{Z}_s - B_3(s)Z_1\rangle\, ds$$
(5.23)
$$\left.\left. - \tfrac{1}{2}\int_0^1 |G(Z_s^x)|^2\, ds\right)\bigg| Z_1 = y - S_1 x\right)$$
$$= \mathbb{E}\exp\left(\sum_{i=0}^{k-1}\langle G(\widehat{Z}_{t_i}^{x,y}), \zeta_{t_{i+1}} - \zeta_{t_i}\rangle\right.$$
$$- \int_{t_i}^{t_{i+1}} \langle G(\widehat{Z}_{t_i}^{x,y}), B_1(s)\widehat{Z}_s - B_3(s)(y - S_1 x)\rangle\, ds$$
$$\left. - \tfrac{1}{2}\int_0^1 |G(\widehat{Z}_s^{x,y})|^2\, ds\right)$$
$$=: \mathbb{E}\Phi_k(x,y),$$

since the processes $(\widehat{Z}_t)$ and $(\xi_t)$ are independent of $Z_1$. By (5.18), the dominated convergence theorem yields

$$\Phi_k(x,y) \to \exp\left(\int_0^1 \langle G(\widehat{Z}_s^{x,y}), d\zeta_s\rangle - \tfrac{1}{2}\int_0^1 |G(\widehat{Z}_s^{x,y})|^2\, ds\right.$$
(5.24)
$$\left. - \int_0^1 \langle G(\widehat{Z}_s^{x,y}), B_1(s)\widehat{Z}_s + B_2(s)x - B_3(s)y\rangle\, ds\right)$$
$$=: \Phi(x,y), \qquad \mathbb{P}\text{-a.s.}$$

(possibly, for a subsequence), since $B_2(\cdot)x, B_3(\cdot)y \in L^1(0,1,H)$ for $x \in H$ and $y \in \mathcal{M}$. Proposition 4.9 and Gaussianity of the process $\widehat{Z}$ imply

(5.25)
$$\mathbb{E}\exp\left(M\int_0^1 |B_1(s)\widehat{Z}_s|\, ds\right) < \infty$$

for each $M < \infty$, and hence, the sequence $(\Phi_k(x,y))$ is equiintegrable on $\Omega$. It follows that

(5.26)
$$\lim_{k\to\infty} \Phi_k(x,y) = \Phi(x,y) \qquad \text{in } L^1(\Omega)$$

and in virtue of (5.23) and (5.20), we find that, for any $x \in E$,

(5.27)
$$\frac{dP(1,x,\cdot)}{d\mu_1^x}(y) = \mathbb{E}\Phi(x,y) = h(x,y), \qquad \mu_1\text{-a.e.}$$

Now we drop the assumption of boundedness of $G$. Proceeding as above, we obtain, for each $N > 0$ and $x \in E$,

$$
\mathbb{E}(\mathbf{1}_{\{\|Z^x\|_{C(0,1;E)} \leq N\}} \exp \rho(Z^x) | Z_1^x = y)
$$
$$
(5.28) \qquad = \mathbb{E}\mathbf{1}_{\{\|\widehat{Z}^{x,y}\|_{C(0,1;E)} \leq N\}} \Phi(x,y), \qquad \mu_1\text{-a.e.},
$$

because on the set $\{\|Z^x\|_{C(0,1;E)} \leq N\}$ we have $G(Z^x) = G_N(Z^x)$ [cf. definition in (5.4)] and $G_N$ is bounded. By (5.3) and the dominated convergence theorem, we get

$$
(5.29) \quad \lim_{N \to \infty} \mathbf{1}_{\{\|Z^x\|_{C(0,1;E)} \leq N\}} \exp \rho(Z^x) = \exp \rho(Z^x) \qquad \text{in } L^1(\Omega),
$$

and (possibly, for a subsequence) it follows that

$$
\lim_{N \to \infty} \mathbb{E}(\mathbf{1}_{\{\|Z^x\|_{C(0,1;E)} \leq N\}} \exp(\rho(Z^x)) | Z_1^x = y)
$$
$$
(5.30) \qquad = \mathbb{E}(\exp(\rho(Z^x)) | Z_1^x = y) = \frac{dP(1,x,\cdot)}{d\mu_1^x}(y), \qquad \mu_1\text{-a.s.}
$$

On the other hand, by the monotone convergence theorem,

$$
(5.31) \quad \lim_{N \to \infty} \mathbb{E}\mathbf{1}_{\{\|\widehat{Z}^{x,y}\|_{C(0,1;E)} \leq N\}} \Phi(x,y) = \mathbb{E}\Phi(x,y) = h(x,y)
$$

for $x \in E$, $y \in \mathcal{M}$, so (5.28), (5.30) and (5.31) yield (5.16) and the proof is completed. □

By means of the formula (5.15), we may find a useful lower estimate on transition densities.

THEOREM 5.3. *For $x \in E$,*

$$
(5.32) \quad \frac{dP(1,x,\cdot)}{d\mu_1}(y) \geq c_1 \exp(-c_2|x|_E^p - \Lambda(y)), \qquad \mu_1\text{-}a.e.,
$$

*where $\Lambda : \mathcal{M}_1 \to \mathbb{R}_+$ is a measurable mapping, $\mathcal{M}_1 \in \mathcal{B}(E), \mu_1(\mathcal{M}_1) = 1$, $p = \max(2, 2m)$ and the constants $c_1, c_2 > 0$ depend only on $A, Q$ and $K, m$ from Hypothesis 2.4(b).*

PROOF. From (5.15) in virtue of the Jensen inequality, we obtain, for $x \in E$,

$$
\frac{dP(1,x,\cdot)}{d\mu_1}(y)
$$
$$
(5.33) \quad \geq \exp\bigg(\mathbb{E}\bigg(\rho(\widehat{Z}^{x,y}) - \int_0^1 \langle G(\widehat{Z}_s^{x,y}), B_1(s)\widehat{Z}_s + B_2(s)x - B_3(s)y \rangle \, ds
$$
$$
+ \langle x, S_1^* Q_1^{-1} y \rangle - \frac{1}{2} |Q_1^{-1/2} S_1 x|^2 \bigg)\bigg),
$$

for $y$ from a set of $\mu_1$-full measure in $E$. Note that in the proof of Theorem 5.2, we found a set $\mathcal{M}, \mu_1(\mathcal{M}) = 1$, such that $\overline{B_3(\cdot)y} \in L^1(0, 1 : H)$ and $h(x, y)$ is well defined for $y \in \mathcal{M}_1$. Similarly, $\overline{S_1^* Q_1^{-1/2}}$ is a Hilbert–Schmidt operator by (2.4), hence, $\overline{S_1^* Q_1^{-1/2}} Q_1^{-1/2} y \in H$ is well defined for $y \in \mathcal{M}_2$, $\mu_1(\mathcal{M}_2) = 1$ and the density $g(x, y)$ is given by the formula (5.12) for $y \in \mathcal{M}_2$. We may take $\mathcal{M}_1 = \mathcal{M} \cap \mathcal{M}_2$. It follows from Hypothesis 2.4(b) and (4.7) that the stochastic integral in $\rho(\widehat{Z}^{x,y})$ is a martingale and, hence, for any $x \in E$,

$$\frac{dP(1, x, \cdot)}{d\mu_1}(y)$$
$$\geq \exp\biggl(-\frac{1}{2} \int_0^1 \mathbb{E}|G(\widehat{Z}_s^{x,y})|^2 \, ds$$
$$- \mathbb{E} \int_0^1 |G(\widehat{Z}_s^{x,y})|(|B_1(s)\widehat{Z}_s^{x,y}| + |B_2(s)x| + |B_3(s)y|) \, ds$$
$$- |x|\overline{|S_1^* Q_1^{-1/2}} Q_1^{-1/2} y| - \frac{1}{2}|Q_1^{-1/2} S_1 x|^2\biggr)$$
$$\geq \exp\biggl(-K^2\biggl(1 + \int_0^1 \mathbb{E}|\widehat{Z}_s^{x,y}|_E^{2m} \, ds\biggr)$$
$$- \mathbb{E} \int_0^1 K(1 + |\widehat{Z}_s^{x,y}|_E^m)(|B_1(s)\widehat{Z}_s| + |B_2(s)x| + |B_3(s)y|) \, ds$$
$$- \tilde{c}|x|_E |S_1^* Q_1^{-1} y| - \frac{1}{2}\tilde{c}^2 \|Q_1^{-1/2} S_1\|^2 \cdot |x|_E^2\biggr),$$
$$x \in E,$$

for $\mu_1$-almost all $y \in \mathcal{M}_1$, where $\tilde{c}$ is the constant from continuous embedding $E \hookrightarrow H$. Set $U_1, U_2 : \mathcal{M}_1 \to \mathbb{R}_+$, $U_1(y) := \|B_3(\cdot)y\|_{L^1(0,1:H)}, U_2(y) = |S_1^* Q_1^{-1} y|$; by (4.7) of Theorem 4.3, we further get

$$\frac{dP(1, x, \cdot)}{d\mu_1}(y) \geq \exp\biggl(-K^2[1 + L(2m)(1 + |x|_E^{2m} + (U(y))^{2m})]$$
$$- K\mathbb{E}\int_0^1 |B_1(s)\widehat{Z}_s| \, ds$$
$$- K\int_0^1 |B_2(s)x| \, ds - KU_1(y)$$
$$- KL(m)(1 + |x|_E^m + U(y)^m)\int_0^1 |B_2(s)x| \, ds$$
(5.34)
$$- KL(m)(1 + \|x\|^m + U^m(y))U_1(y)$$

$$- KL(m)\mathbb{E}\|\widehat{Z}^{x,y}\|_{C(0,1;E)}^m \int_0^1 |B_1(s)\widehat{Z}_s|\, ds$$

$$- \tilde{c}|x|_E U_2(y) - \frac{1}{2}\tilde{c}^2\|Q_1^{-1/2}S_1\|^2 \cdot |x|_E^2\biggr).$$

By Lemma 3.2, we have

$$\begin{aligned}
\int_0^1 |B_2(s)x|\, ds &= \|BQ_1^{-1/2}S_1 x\|_{L^1(0,1;H)} \\
&\leq \|BQ_1^{-1/2}S_1 x\|_{L^2(0,1;H)} \\
&= |Q_1^{-1/2}S_1 x| \leq \tilde{c}\|Q_1^{-1/2}S_1\| \cdot |x|_E
\end{aligned} \tag{5.35}$$

and it follows from Proposition 4.9 that

$$\mathbb{E}\left(\int_0^1 |B_1(s)\widehat{Z}_s|\, ds\right)^q < \infty, \tag{5.36}$$

for each $q < \infty$. Therefore, for each $\eta > 0$ small enough, there exist constants $c_1(\eta) > 0$ and $c_3(\eta) > 0$ and a function $\Lambda = \Lambda_\eta : \mathcal{M}_1 \to \mathbb{R}_+$ such that, $x \in E$ and $y \in \mathcal{M}_1$,

$$\begin{aligned}
\frac{dP(1,x,\cdot)}{d\mu_1}(y) &\geq \exp\Bigl(-c_1(\eta) - (K^2 L(2m) + \eta)|x|_E^{2m} \\
&\quad - KL(m)\|Q_1^{-1/2}S_1\| \cdot |x|_E^{m+1} \\
&\quad - c_3(\eta)|x|_E^{m+\eta} - \left(\frac{1}{2}\tilde{c}^2\|Q_1^{-1/2}S_1\|^2 + \eta\right)|x|_E^2 - \Lambda(y)\Bigr),
\end{aligned}$$

and the estimate (5.32) follows. $\square$

REMARK 5.4. Under more stringent conditions, we may obtain a lower estimate on the transition density which is more "explicit" in $y$ and has a more symmetric form. In addition to the conditions of Theorem 5.2, assume that there exists a Banach space $\tilde{E}$ of $\mu_1$-full measure, continuously embedded into $H$ such that

$$S_1(\tilde{E}) \subset \mathrm{im}(Q_1), \tag{5.37}$$

$$\int_0^1 |B_3(s)y|\, ds \leq a_1 |y|_{\tilde{E}}, \qquad y \in \tilde{E}, \tag{5.38}$$

and

$$U(y) = \sup_{t\in[0,1]}\left|\int_0^t S_{t-s}QS_{1-s}^* Q_1^{-1} y\, ds\right|_E \leq a_2 |y|_{\tilde{E}}, \qquad y \in \tilde{E}, \tag{5.39}$$

for some $a_1, a_2 > 0$. Then for some constants $b_1, b_2, b_3 > 0$ (dependent only on $A, Q, K$ and $m$) and $p = \max(2, 2m)$, we have

$$(5.40) \quad \frac{dP(1, x, \cdot)}{d\mu_1}(y) \geq b_1 \exp\{-b_2 |x|_E^p - b_3 |y|_{\tilde{E}}^p\}, \qquad x \in E, y \in \tilde{E} \text{ a.e.}$$

To see (5.40), we check that, under present conditions (5.37)–(5.39), (5.34) implies, for all for $x \in E$ and for $\mu_1$-a.e. $y \in \tilde{E}$,

$$\frac{dP(1, x, \cdot)}{d\mu_1}(y) \geq \exp(-C(1 + |x|_E^{2m} + |y|_{\tilde{E}}^{2m} - |x|_E + |y|_{\tilde{E}}$$
$$+ |x|_E^{m+1} + |y|_{\tilde{E}}^m |x|_E + |x|_E^m |y|_{\tilde{E}}$$
$$+ |x|_E^{m+\eta} + |y|_{\tilde{E}}^{m+\eta} + |x|_E^2 + |y|_{\tilde{E}}^2)),$$

for arbitrary small $\eta > 0$ and a universal $C = C(\eta) < \infty$, and (5.40) follows. It may be of interest to mention some particular cases when the conditions (5.37)–(5.39) are satisfied. A trivial example is a finite-dimensional one, $H = E = \tilde{E} = \mathbb{R}^d$, in which case we obtain

$$(5.41) \quad \frac{dP(1, x, \cdot)}{d\mu_1}(y) \geq b_1 \exp(-b_2 |x|_{\mathbb{R}^d}^p - b_3 |y|_{\mathbb{R}^d}^p), \qquad x, y \in \mathbb{R}^d.$$

Note that in this case the only assumptions in Theorem 5.2, Corollary 5.3 and the present remark are the well posedness and growth conditions in Hypothesis 2.4 and the strong Feller property of the Ornstein–Uhlenbeck process (2.7).

Suppose that $A = A^*$ is strictly negative and define $H_\lambda = \text{dom}((-A)^\lambda), \lambda \geq 0$, with the norm $|y|_\lambda = |(-A)^\lambda y|, y \in \text{dom}((-A)^\lambda)$. Let $Q = I$; then $A^{-1}$ must be compact and it is easy to check that $\text{im}(S_1) \subset \text{im}(Q_1) = \text{dom}(A)$, therefore, (5.37) holds with any $\tilde{E}, \tilde{E} \hookrightarrow H$. Furthermore, we have

$$(5.42) \quad \begin{aligned} |Q^{1/2} S_{1-s} Q_1^{-1} y| &= \|S_{1-s}(-A)^{1-\lambda}\| \|(-A)^{\lambda-1} Q_1^{-1} (-A)^{-\lambda}\| \cdot |y|_\lambda \\ &\leq \frac{\text{const}}{(1-s)^{1-\lambda}} \|y\|_\lambda \end{aligned}$$

for $y \in H_\lambda$ since $A^{-1} Q_1^{-1} \in \mathcal{L}(H)$, thus, (5.38) holds for $\tilde{E} = H_\lambda$ with any $\lambda > 0$. If, in addition, $\|S_t\|_{\mathcal{L}(H, E)} \leq \text{const} \cdot t^{-\sigma}, t \in [0, 1]$, for $\sigma > 0$ such that $\sigma < \lambda$, then (5.39) holds as well since

$$(5.43) \quad |S_{t-s} S_{1-s} Q_1^{-1} y|_E \leq \text{const}(t-s)^{-\sigma}(1-s)^{\lambda-1}, \qquad 0 < s < t \leq 1.$$

**6. Exponential convergence to invariant measure.** The following uniform ultimate moment boundedness result will be useful in the sequel.

PROPOSITION 6.1. *Assume that the growth condition* (2.11) *holds true and*

(6.1) $$k(p) := \sup_{t \geq 0} \mathbb{E}|Z_t|_E^p < \infty, \qquad p > 0.$$

*Then*

(6.2) $$\mathbb{E}_x|X_t|_E \leq e^{-k_1 t}|x|_E + \frac{k_2 k(s) + k_3}{k_1} + k(1), \qquad t \geq 0.$$

*Suppose that the following stronger version of* (2.11) *holds*: *For each* $x \in \mathrm{dom}(\tilde{A})$, *there exists* $x^* \in \partial |x|_E$ *such that, for some* $k_1, k_2, k_3 > 0, s > 0, \varepsilon > 0$, *we have*

(6.3) $\langle \tilde{A}x + F(x+y), x^* \rangle_{E,E^*} \leq -k_1 |x|_E^{1+\varepsilon} + k_2 |y|_E^s + k_3, \qquad y \in E.$

*Then*

(6.4) $$\sup_{x \in E} \sup_{t \geq 1} \mathbb{E}_x |X_t|_E \leq \widehat{M},$$

*where*

(6.5) $$\widehat{M} = k(1) + \max\left(\left(\frac{2(k_2 k(s) + k_3)}{k_1}\right)^{1+\varepsilon}, \left(\frac{1}{k_1 \varepsilon} + 2\right)^{1/\varepsilon}\right).$$

PROOF. Inequality (6.4) has been proven in [14], Proposition 2.1 (see also a similar result in [20]). The proof of (6.2) follows the lines of similar proofs based on Yosida approximation techniques (see, e.g., [9]) and we sketch it only. The process $Y^x(t) := X_t^x - Z_t$ satisfies the equation

(6.6) $$Y^x(t) = S_t x + \int_0^t S_{t-s} F(Y^x(s) + Z_s)\, ds, \qquad t \geq 0,$$

and the sequence of approximating processes $Y_\lambda(t)$ is defined by

(6.7) $$Y_\lambda^x(t) = R(\lambda) S_t x + \int_0^t R(\lambda) S_{t-s} F(Y^x(s) + Z_s)\, ds, \qquad t \geq 0,$$

where $R(\lambda) := \lambda(\lambda I - \tilde{A})^{-1} \in \mathcal{L}(E)$ is well defined for $\lambda$ large enough. It is well known that

(6.8) $Y_\lambda^x \to Y^x, \qquad \dfrac{dY_\lambda}{dt} - \tilde{A} Y_\lambda - F(Y_\lambda + Z) = \sigma_\lambda^x \to 0, \qquad \lambda \to \infty,$

in $C(0, T; E)$ (cf. page 201 of [9]). Since by (2.11)

(6.9) $$\frac{d^-}{dt}|Y_\lambda^x(t)|_E \leq -k_1 |Y_\lambda(t)|_E + k_2 |Z_t^x|_E^s + k_3 + |\sigma_\lambda(t)|_E,$$

we obtain

$$|Y^x(t)|_E \leq e^{-k_1 t}|x|_E + \int_0^t e^{-k_1(t-\tau)}(k_2 |Z_\tau^x|_E^s + k_3)\, d\tau,$$

and thereby,

$$\mathbb{E}|Y^x(t)|_E \leq e^{-k_1 t}|x|_E + \int_0^t e^{k_1(t-\tau)}(k_2 k(s) + k_3)\, d\tau \tag{6.10}$$

and (6.2) follows. $\square$

Our next aim is to establish uniform geometric ergodicity and $V$-uniform ergodicity results for $V(x) = |x|_E + 1$ using the lower density estimates and uniform moment boundedness shown above. We will also find explicit bounds on the convergence rates, hence, the constants below will play some role. We assume that

$$\mathbb{E}_x|X_t|_E \leq k_0 e^{-k_1 t}|x|_E + \widehat{c}, \qquad t \geq 0, \tag{6.11}$$

for some $k_0$, $k_1 > 0$ and $\widehat{c} \in \mathbb{R}$. Note that, by Proposition 6.1, if (6.1) and the growth condition (2.11) are both satisfied, then (6.11) holds with $k_1$ given in (2.11), $k_0 = 1$ and

$$\widehat{c} = \frac{k_2 k(s) + k_3}{k_1} + k(1). \tag{6.12}$$

Now, take $R > 4\widehat{c}$, $r > 4(\widehat{c} + \frac{1}{2})$, and define

$$\begin{aligned}
t_0 &= -\frac{1}{k_1}\log\left(\frac{R}{2rk_0} - \frac{\widehat{c}}{rk_0}\right), \\
T &= \max\left(t_0 + 1, -\frac{1}{k_1}\log\frac{1}{4k_0}\right), \qquad b = \widehat{c} + \frac{1}{2}
\end{aligned} \tag{6.13}$$

and

$$\delta = \tfrac{1}{2} c_1 e^{-c_2 R^p} \int_{B_r} e^{-\Lambda(y)} \mu_1(dy), \tag{6.14}$$

where $B_r := \{y \in E, |y|_E < r\}$, and $c_1, c_2, p$ and $\Lambda$ are defined in the same way as in Theorem 5.3. In the following proposition, existence of a universal small set satisfying a uniform geometric drift condition is shown.

PROPOSITION 6.2. *Assume* (6.11). *Then the following holds:*

(a) *We have*

$$\inf_{x \in B_r} P(T, x, \Gamma) \geq \delta \bar{\mu}(\Gamma), \qquad \Gamma \in \mathcal{B}(E), \tag{6.15}$$

*where*

$$\bar{\mu}(\Gamma) := \left(\int_{B_r} e^{-\Lambda(y)} \mu_1(dy)\right)^{-1} \int_{B_r \cap \Gamma} e^{-\Lambda(y)} \mu_1(dy), \qquad \Gamma \in \mathcal{B}(E), \tag{6.16}$$

*is a probability measure. In particular, $B_r$ is a small set of the Markov chain $(\widetilde{X}_n) := (X_{nT})$, with the lower bound measure $\delta\bar{\mu}$.*

(b) *We have*

(6.17) $$\mathbb{E}_x(|X_T|_E + 1) \leq \tfrac{1}{2}(|x|_E + 1) + b\mathbf{1}_{B_r}(x), \qquad x \in E,$$

*that is, the chain* $(\tilde{X}_n)$ *satisfies the one-step Lyapunov–Foster condition of geometric drift toward* $B_r$, *with the constants* $\tfrac{1}{2}$ *and* $b$ *and the Lyapunov function* $V(x) = |x|_E + 1$.

PROOF. (a) By (6.11), we have, for $t \geq t_0, x \in E, |x|_E \leq r$,

(6.18) $$P(t, x, B_R) \geq 1 - \frac{\mathbb{E}|X_t|_E}{R} \geq 1 - \frac{1}{R}(k_0 r e^{-k_1 t} + \widehat{c}) \geq \frac{1}{2}$$

and therefore, by Theorem 5.3, for each $t \geq t_0$, we get

(6.19) $$\begin{aligned} P(t+1, x, \Gamma) &= \int_E P(1, y, \Gamma) P(t, x, dy) \geq \int_{B_R} P(1, y, \Gamma) P(t, x, dy) \\ &\geq \int_{B_R} \int_\Gamma \frac{dP(1, y, \cdot)}{d\mu_1}(z) \mu_1(dz) P(t, x, dy) \\ &\geq \int_{B_R} \int_\Gamma c_1 \exp\{-c_2 |y|_E^p - \Lambda(z)\} \mu_1(dz) P(t, x, dy) \\ &\geq c_1 e^{-c_2 R^p} \int_\Gamma e^{-\Lambda(z)} \mu_1(dz) P(t, x, B_R), \qquad x \in E, \Gamma \in \mathcal{B}(E). \end{aligned}$$

Hence,

(6.20) $$\begin{aligned} \inf_{x \in B_r} P(t+1, x, \Gamma) &\geq c_1 e^{-c_2 R^p} \int_\Gamma e^{-\Lambda(z)} \mu_1(dz) \inf_{x \in B_r} P(t, x, B_R) \\ &\geq \tfrac{1}{2} c_1 e^{-c_2 R^p} \int_{\Gamma \cap B_r} e^{-\Lambda(z)} \mu_1(dz) \\ &= \delta \bar{\mu}(\Gamma), \qquad \Gamma \in \mathcal{B}(E) \end{aligned}$$

and (6.15) follows.

To prove part (b), we use again (6.11) to obtain

(6.21) $$\begin{aligned} \mathbb{E}_x(|X_t|_E + 1) &\leq k_0 |x|_E e^{-k_1 t} + \widehat{c} + 1 \leq \tfrac{1}{4} |x|_E + \widehat{c} + 1 \\ &\leq \tfrac{1}{2}(|x|_E + 1) - \tfrac{1}{4}|x|_E - \tfrac{1}{2} + \widehat{c} + 1 \\ &\leq \tfrac{1}{2}(|x|_E + 1) + (\widehat{c} + \tfrac{1}{2}) \mathbf{1}_{B_r}(x) \end{aligned}$$

for $x \in E$, $t \geq -\tfrac{1}{k_1} \log \tfrac{1}{4k_0}$, which completes the proof. □

In the next theorem our main result on uniform geometric $V$-ergodicity for $V(x) = |x|_E + 1$ is stated. It is based on the paper by Meyn and Tweedie [33], where exact bounds for geometric ergodicity of irreducible Markov chains are

found, and Proposition 6.2 above. Following [33], we introduce the constants $v, M_c, \gamma_c, \widehat{\lambda}, \widehat{b}$ and $\bar{\xi}$ as follows:

(6.22)
$$v = r+1, \qquad \gamma_c = \delta^{-2}(4b + \delta v), \qquad \widehat{\lambda} = \frac{1/2 + \gamma_c}{1+\gamma_c} < 1,$$
$$\widehat{b} = v + \gamma_c, \qquad \bar{\xi} = \frac{4-\delta^2}{\delta^5} 4b^2$$

and

$$M_c = \frac{1}{(1-\widehat{\lambda})^2}(1 - \widehat{\lambda} + \widehat{b} + \widehat{b}^2 + \bar{\xi}(\widehat{b}(1-\widehat{\lambda}) + \widehat{b}^2)) > 1.$$

We will show that the Markov chain $(\tilde{X}_n)$ has the geometric rate of convergence to the invariant measure with any constant

(6.23)
$$\rho \in \left(1 - \frac{1}{M_c}, 1\right).$$

Let $b_V \mathcal{B}$ denote the Banach space of measurable functions $\varphi : E \to \mathbb{R}$ such that

$$\|\phi\|_V = \sup_{x \in E} \frac{|\varphi(x)|}{V(x)} < \infty.$$

THEOREM 6.3. *Assume* (6.11). *Then there exists an invariant measure* $\mu^* \in \mathcal{P}$ *and for* $V(x) = |x|_E + 1$, *we have*

(6.24)
$$\sup_{\|\phi\|_V \leq 1} \left| P_t \varphi(x) - \int_E \varphi \, d\mu^* \right| \leq MV(x) e^{-\omega t}, \qquad t \geq 0, x \in E,$$

*where*

(6.25)
$$\omega = -\frac{1}{T} \log \rho > 0 \quad and$$
$$M = (1+\gamma_c) \frac{\rho}{\rho + M_c^{-1} - 1} (\widehat{c} + k_0 + 1) e^{-\log \rho}$$

*and*

(6.26)
$$\|P_t^* \nu - \mu^*\|_{\mathrm{var}} \leq M(L_\nu + 1) e^{-\omega t}, \qquad t \geq 0, \nu \in \mathcal{P},$$

*where* $L_\nu = \int_E |x|_E \nu(dx)$. *The constants* $\omega$ *and* $M$ *may be chosen the same for all nonlinear terms* $F$ *satisfying Hypothesis* 2.4(b) *with the same constants* $K$ *and* $m$ *and* (6.11) *with the same constants* $k_0, k_1$ *and* $\widehat{c}$ [*or, in particular, satisfying the growth condition* (2.11) *with the same* $k_1, k_2, k_3$ *and* $s$].

PROOF. The existence of an invariant probability measure in presence of the lower bound measure [cf. (6.15)] and condition (6.11) are well known (see, e.g., [28]). It follows from Proposition 6.2 that we may apply Theorem 2.3 of [33] to the Markov chain $(\tilde{X}_n) = (X_{nT})$ (note that the measure $\bar{\mu}$ is concentrated on $B_r$, so the conditions imposed in [33] are satisfied), which yields

$$\sup_{\|\phi\|_V \leq 1} \left| P_{nT}\varphi(x) - \int \varphi \, d\mu^* \right| \leq C\rho^n V(x) = Ce^{-nT\omega}V(x), \tag{6.27}$$

$$x \in E, n \in \mathbb{N},$$

where $C = (1 + \gamma_c)\frac{\rho}{\rho + M_c^{-1} - 1}$. Using the semigroup property of $(P_t)$ and (6.11), we obtain

$$\begin{aligned}
\sup_{\|\phi\|_V \leq 1} \left| P_{nT+s}\varphi(x) - \int \varphi \, d\mu^* \right| &\leq \sup_{\|\phi\|_V \leq 1} \left| P_s\left(P_{nT}\varphi - \int \varphi \, d\mu^*\right)(x) \right| \\
&\leq |P_s(Ce^{-nT\omega}V(x))| \leq Ce^{-nT\omega}\mathbb{E}_x(|X_s|_E + 1) \\
&\leq Ce^{-nT\omega}(k_0|x|_E + \hat{c} + 1) \\
&\leq C(\hat{c} + k_0 + 1)e^{-\log\rho}(|x|_E + 1)e^{-\omega(nT+s)},
\end{aligned} \tag{6.28}$$

$$n \in \mathbb{N}, x \in E, s \in [0, T],$$

which yields (6.24). Inequality (6.26) is an obvious consequence of (6.24) since

$$\begin{aligned}
\|P_t^*\nu - \mu^*\|_{\mathrm{var}} &\leq \int_E \|P(t, x, \cdot) - \mu^*\|_{\mathrm{var}}\nu(dx) \\
&\leq \int \sup_{\|\phi\|_V \leq 1} \left| P_t\varphi(x) - \int \varphi \, d\mu^* \right| \nu(dx) \\
&\leq \int M(|x|_E + 1)e^{-\omega t}\nu(dx) = M(L_\nu + 1)e^{-\omega t}
\end{aligned} \tag{6.29}$$

for each $\nu \in \mathcal{P}$, $t \geq 0$. The universality of $M$ and $\omega$ follows from the fact that all constants defined in (6.12)–(6.14) and (6.22) (including $c_1, c_2, p$ and the mapping $\Lambda$, cf. Theorem 5.3) are independent of $F$. □

If the growth of the nonlinear term $F$ is faster than linear, the Markov process defined by the equation (2.1) may be uniformly ergodic, that is, the constant $L_\nu$ in (6.26) may be replaced by another constant independent of the initial measure $\nu \in \mathcal{P}$. This has been established earlier in [14] and [20]; however, the lower bound measures are not found there constructively. In the theorem below explicit bounds are found and, in particular, uniformity of convergence with respect to coefficients is proven.

THEOREM 6.4. *Assume* (6.1) *and let the stronger growth condition* (6.3) *hold true. Then there exists an invariant measure $\mu^* \in \mathcal{P}$ and for any $\nu \in \mathcal{P}$,*

$$(6.30) \qquad \|P_t^* \nu - \mu^*\|_{\mathrm{var}} \leq (1-\delta)^{-1} e^{-\widehat{\omega} t} \|\nu - \mu^*\|_{\mathrm{var}}, \qquad t \geq 0,$$

*where $\widehat{\omega} = -\frac{1}{2} \log(1-\delta) > 0$, $\delta$ is defined by* (6.14) *with $R = 2\widehat{M}$ and $r = \infty$ and $\widehat{M}$ is given by* (6.5). *In particular, the constants on the r.h.s. of* (6.30) *are uniform with respect to all nonlinear terms $F$ satisfying the growth conditions* (2.9) *and* (6.3) *with the same constants $K, m, k_1, k_2, k_3, s$ and $\varepsilon$.*

PROOF. By Proposition 6.1, we have that

$$(6.31) \qquad \inf_{x \in E} P(1, x, B_R) \geq 1 - \sup_{x \in E} \frac{\mathbb{E}_x |X_1|_E}{R} \geq \frac{1}{2}$$

and similarly, as in (6.19), we get

$$(6.32) \quad \begin{aligned} \inf_{x \in E} P(2, x, \Gamma) &\geq \inf_{x \in E} \int_E P(1, y, \Gamma) P(1, x, dy) \\ &\geq \tfrac{1}{2} c_1 e^{-c_2 R^p} \int_\Gamma e^{-\Lambda(z)} \mu_1(dz) = \delta \bar{\mu}(\Gamma), \end{aligned}$$

where $\bar{\mu}$ is defined by (6.16) with $r = \infty$. For each $\nu \in \mathcal{P}$, it follows that $P_2^* \nu \geq \delta \bar{\mu}$ and a simple computation (cf., e.g., [14], Theorem 2.4) yields

$$\|P_2^* \mu\|_{\mathrm{var}} \leq (1-\delta) \|\mu\|_{\mathrm{var}}, \qquad \mu = \nu_1 - \nu_2, \nu_1, \nu_2 \in \mathcal{P}.$$

By the semigroup property of $(P_t^*)$, we have $\|P_{2n}^* \mu\|_{\mathrm{var}} \leq (1-\sigma)^n \|\mu\|_{\mathrm{var}}$ and for $s \in [0, 2]$, it follows that

$$\begin{aligned} \|P_{2n+s}^* \mu\|_{\mathrm{var}} &\leq \|P_s^* P_{2n}^* \mu\|_{\mathrm{var}} \leq \|P_{2n}^* \mu\|_{\mathrm{var}} \\ &\leq e^{-2n\widehat{\omega}} \|\mu\|_{\mathrm{var}} \leq (1-\delta)^{-1} e^{-(2n+s)\widehat{\omega}} \|\mu\|_{\mathrm{var}}. \end{aligned} \qquad \square$$

**7. Uniform spectral gap property.** In this section we consider exponential ergodicity in spaces $L^p(E, \mu^*)$ for $p \in [1, \infty)$. Note first that, by Theorem 5.2, the transition kernels $P(T, x, \cdot)$ are equivalent for $T > 0$, $x \in E$, they are also equivalent to the invariant measure $\mu^*$ (if it exists) and we have, for each $t > 0$,

$$P_t \phi(x) = \int_E p_t(x, y) \phi(y) \mu^*(dy), \qquad \phi \in C_b(E),$$

where the function $(x, y) \to p_t(x, y)$ is measurable. Let us recall that $(P_t)$ extends to a contraction semigroup on $L^p(E, \mu^*)$ for all $p \in [1, \infty]$ and is a $C_0$-semigroup if $p < \infty$. Let $\mathcal{M}(E) \subset (C_b(E))^*$ denote the space of finite Borel measures on $E$ with the variation norm. For $\nu \in \mathcal{M}(E)$, we have

$$(7.1) \qquad \langle P_t^* \nu, \phi \rangle = \langle \nu, P_t \phi \rangle = \int_E \int_E p_t(x, y) \phi(y) \mu^*(dy) \nu(dx).$$

LEMMA 7.1. *Assume that the equation* (2.1) *has an invariant measure* $\mu^* \in \mathcal{P}$. *Then the space* $L^p(E, \mu^*)$ *is invariant for* $(P_t^*)$ *for each* $p \in [1, \infty)$. *Moreover,* $\|P_t^*\|_{p \to p} = 1$ *and*

$$P_t^* \psi(y) = \int_E p_t(x, y) \psi(x) \mu^*(dx), \qquad \psi \in L^1(E, \mu^*). \tag{7.2}$$

*Finally,* $(P_t^*)$ *is a* $C_0$-*semigroup on* $L^p(E, \mu^*)$ *for* $p \in [1, \infty)$.

PROOF. For $\psi \in L^1(E, \mu^*)$, $\psi \geq 0$, we define

$$G_t \psi(y) = \int_E p_t(x, y) \psi(x) \mu^*(dx),$$

and $\nu = \psi \mu^*$. For $\phi \geq 0$, the Fubini theorem yields

$$\langle P_t^* \nu, \phi \rangle = \int_E \int_E p_t(x, y) \psi(x) \phi(y) \mu^*(dx) \mu^*(dy) = \langle G_t \psi, \phi \rangle < \infty.$$

Putting $\phi = 1$, we obtain

$$\|G_t \psi\|_1 = \langle P_t^*(\psi \mu^*), 1 \rangle = \langle \psi \mu^*, 1 \rangle = \|\psi\|_1,$$

and therefore,

$$\|G_t \psi\|_1 \leq \|\psi\|_1, \qquad \psi \in L^1(E, \mu^*), \psi \geq 0.$$

Clearly, $G_t \psi = P_t^*(\psi \mu^*)$. All those arguments extend immediately to an arbitrary $\psi \in L^1(E, \mu^*)$ and therefore, $G_t$ is a contraction on $L^1(E, \mu^*)$. Other parts follow easily by a standard density argument. $\square$

Let $L_p$ be the generator of $(P_t)$ acting in $L^p(E, \mu^*)$. We say that $L_p$ has the spectral gap in $L^p(E, \mu^*)$ if there exists $\delta > 0$ such that

$$\sigma(L_p) \cap \{\lambda : \operatorname{Re} \lambda > -\delta\} = \{0\}.$$

The largest $\delta$ with this property will be denoted by $\operatorname{gap}(L_p)$.

THEOREM 7.2. *Assume* (6.1) *and let the stronger growth condition* (6.3) *be satisfied. Then, for each* $p \in (1, \infty)$, *we have*

$$\operatorname{gap}(L_p) \geq \frac{\widehat{\omega}}{p} \tag{7.3}$$

*and*

$$\|P_t \phi - \langle \mu^*, \phi \rangle\|_p \leq C_p e^{-(\widehat{\omega}/p)t} \|\phi\|_p, \tag{7.4}$$

*where* $\widehat{\omega} > 0$ *is defined in Theorem* 6.4. *If, moreover, the semigroup* $(P_t)$ *is symmetric in* $L^2(E, \mu^*)$, *then* (7.3) *and* (7.4) *hold for* $p = 1$.

PROOF. By Theorem 6.4, there exist $C > 0$ such that

$$\|P_t^* \nu - \mu^*\|_{\text{var}} \leq Ce^{-\widehat{\omega}t}, \tag{7.5}$$

for any probability measure $\nu$ on $E$ or, equivalently,

$$\|P_t^* \nu - \nu(E)\mu^*\|_{\text{var}} \leq C\|\nu\|_{\text{var}} e^{-\widehat{\omega}t}, \tag{7.6}$$

for any signed measure $\nu$. If $\nu = \psi\mu^*$, then (7.6) and Lemma 7.1 imply

$$\|P_t^* \psi - \langle \psi, 1 \rangle\|_1 \leq C\|\psi\|_1 e^{-\widehat{\omega}t},$$

hence,

$$\|P_t^* - \Pi\|_1 \leq Ce^{-\widehat{\omega}t}, \tag{7.7}$$

where $\Pi\psi = \langle \psi, 1 \rangle 1$ and

$$\|P_t^* - \Pi\|_{q \to q} \leq 2, \qquad q \in [1, \infty]. \tag{7.8}$$

Take $q > 2$. Then by (7.7), (7.8) and the Riesz–Thorin theorem, we find that

$$\|P_t^* - \Pi\|_{p \to p} \leq C^\theta e^{-\widehat{\omega}\theta t} 2^{1-\theta}, \tag{7.9}$$

where

$$\frac{1}{p} = \frac{\theta}{1} + \frac{1-\theta}{q}.$$

Therefore, taking $q \to \infty$ in (7.9), we obtain

$$\|P_t^* - \Pi\|_{p \to p} \leq C_2 e^{-(\widehat{\omega}/p t)}. \tag{7.10}$$

Therefore, (7.4) holds, and since $(P_t)$ is a $C_0$-semigroup in $L^2(E, \mu^*)$, Theorem 3.6.2 in [34] (7.10) implies the spectral gap property with $\text{gap}(L_p) \geq \frac{\widehat{\omega}}{p}$ for $p \in (1, \infty)$. If $(P_t)$ is symmetric, then the conclusion of the theorem for $p = 1$ follows immediately from (7.5). $\square$

REMARK 7.3. (1) In Theorem 7.1 and 7.2 the invariant measure $\mu^*$, hence, the space $L^p(E, \mu^*)$, depends on the coefficients of equation (2.1). It is interesting to note that the lower bound on $\text{gap}(L_p)$ and $C_p$ are universal for all systems satisfying Hypothesis 2.4(b) and (6.3) with the same constants.

(2) By Theorem 7.2, the spectral gap exists for all $p \in (1, \infty)$. The fact that this property holds in $L^1(E, \mu^*)$ is perhaps surprising. Note that it does not need to hold in general if $F = 0$. It is known (cf. [13]) that, for a one-dimensional Ornstein–Uhlenbeck operator $L_1^{\text{OU}}$ considered in $L^1(E, \mu^*)$, we have

$$\sigma(L_1^{\text{OU}}) = \{\lambda : \text{Re}\,\lambda \leq 0\}.$$

If $p=2$ and $(P_t)$ is symmetric, the stronger growth condition (6.3) is not needed. We may get an estimate on spectral gap in $L^2(E,\mu^*)$ under the standard ultimate boundedness assumption, which is stated in Corollary 7.4 below. Note that the assumption of symmetricity of $(P_t)$ may not be removed (cf. Example 9.1 below).

COROLLARY 7.4. *Let the conditions of Theorem 6.3 be satisfied and assume that $(P_t)$ is symmetric on $L^2 = L^2(E,\mu^*)$. Then*

$$(7.11) \qquad \|P_t\varphi\|_{L^2} \leq e^{-\omega t}\|\varphi\|_{L^2}$$

*holds for all $t \geq 0$ and $\varphi \in L^2, \int \varphi\,d\mu^* = 0$, where $\omega$ is defined in (6.25).*

PROOF. Taking arbitrary $t > 0$, we get, by Theorem 6.3,

$$(7.12) \qquad \sup_{\|\phi\|_V \leq 1} \left| P_{tn}\varphi(x) - \int \varphi\,d\mu^* \right| \leq MV(x)\widehat{\rho}^n, \qquad x \in E,$$

where $\widehat{\rho} = e^{-\omega t}$, that is, the skeleton $(Y_n) := (X_{tn})$ is $V$-uniformly ergodic with the rate $\widehat{\rho}$. By [36], Theorem 2.1, it follows that

$$(7.13) \qquad \|P_{tn}\varphi\|_{L^2} \leq \widehat{\rho}^n \|\varphi\|_{L^2}, \qquad n \in \mathbb{N}, \varphi \in L^2, \int \varphi\,d\mu^* = 0$$

and taking $n = 1$, we obtain (7.11). □

## 8. Some extensions.

8.1. *Equations nonhomogeneous in time.* Some results in the present paper may be easily generalized to the case when the nonlinear term $F = F(t,x)$ in the equation (2.1) also depends on time, that is, the equation has the form

$$dX_t = (AX_t + F(t,X_t))\,dt + \sqrt{Q}\,dW_t, \qquad t \geq s \geq 0,$$
$$X_s = x,$$

and defines a nonhomogeneous Markov process. For instance, let $P_{s,t}$ and $P^*_{s,t}$ denote the corresponding two-parameter Markov semigroup and adjoint Markov semigroup, respectively, $0 \leq s \leq t$, and set $P(s,x,t,\Gamma) := \mathbb{E}_{s,x}\mathbf{1}_\Gamma(X_t) = P_{s,t}\mathbf{1}_\Gamma(x), 0 \leq s \leq t, x \in E, \Gamma \in \mathcal{B}$.

THEOREM 8.1. *Let Hypotheses 2.1, 2.2, 2.3 and condition (6.1) be satisfied and let $F: \mathbb{R}_+ \times E \to E$ be a jointly measurable mapping such that $F(t,\cdot)$ is Lipschitz continuous on bounded sets in $E$ and satisfies Hypothesis 2.4(b) and the growth condition (6.3) with constants independent of $t \in \mathbb{R}_+$. Then*

$$\|P^*_{s,t}\nu_1 - P^*_{s,t}\nu_2\|_{\text{var}} \leq (1-\delta)^{-1} e^{-\widehat{\omega}(t-s)} \|\nu_1 - \nu_2\|_{\text{var}}, \qquad 0 \leq s \leq t, \nu_1,\nu_2 \in \mathcal{P},$$

*where $\widehat{\omega} = -\frac{1}{2}\log(1-\delta)$ and $\delta, \widehat{\omega}$ depend only on the constants in Hypothesis 2.4(b) and (6.3).*

The proof is just a a slight modification of the above results; similarly to Theorems 5.2 and 5.3 and Proposition 6.1, we obtain

$$\inf_{s\in\mathbb{R}_+, x\in E} P(s,x,s+2,\Gamma) \geq \inf_{s\in\mathbb{R}_+, x\in E} \int_E P(s+1,y,s+2,\Gamma)P(s,x,s+1,dy)$$

$$\geq \tfrac{1}{2}c_1 e^{-c_2 R^p}\int_\Gamma e^{-\Lambda(z)}\mu_1(dz) = \delta\bar{\mu}(\Gamma),$$

where $R = 2\widehat{M}$ and $\widehat{M}$ is defined in Proposition 6.1. Our statement now follows just as in the proof of Theorem 6.4.

8.2. *Continuous dependence of invariant measures on parameter.* Uniformity of convergences proven in Theorems 6.3 and 6.4 with respect to nonlinear drifts in a fairly large class may be useful in some cases, for instance, in ergodic control theory. Another application (given below) is a continuous dependence of invariant measures on a parameter. Consider the parameter-dependent equation

(8.1)
$$dX_t^\alpha = (AX_t^\alpha + F_\alpha(X_t^\alpha))\,dt + \sqrt{Q}\,dW_t,$$
$$X_0^\alpha = x \in E,$$

where $\alpha \in \mathcal{A} \subset \mathbb{R}^d$. Denote by $P^\alpha(t,x,\cdot)$ and $\mu^\alpha$ the transition probability kernel and the invariant measure, respectively, associated with the equation (8.1).

THEOREM 8.2. *Let Hypotheses* 2.1–2.4, *condition* (6.1) *and the growth condition* (2.11) *hold for equation* (8.1) *with the constants independent of* $\alpha \in \mathcal{A}$, *and assume*

(8.2)
$$\lim_{\alpha\to\alpha_0} G_\alpha(x) = G_{\alpha_0}(x), \qquad x \in E,$$

*where* $F_\alpha = Q^{1/2}G_\alpha$. *Then*

(8.3)
$$\lim_{\alpha\to\alpha_0}\|\mu_\alpha - \mu_{\alpha_0}\|_{\mathrm{var}} = 0.$$

PROOF. First we prove

(8.4)
$$\lim_{\alpha\to\alpha_0}\|P^\alpha(t,x,\cdot) - P^{\alpha_0}(t,x,\cdot)\|_{\mathrm{var}} = 0,$$

for each $t > 0$ and $x \in E$. By (5.2), it suffices to show that $\exp\rho^\alpha \to \exp\rho^{\alpha_0}$ in $L^1(\Omega)$, where

$$\rho_\alpha(Z^x) = \int_0^t \langle G_\alpha(Z_s^x), dW_s\rangle - \tfrac{1}{2}\int_s^t |G_\alpha(Z_s^x)|^2\,ds.$$

By the dominated convergence theorem, we have

$$\lim_{\alpha \to \alpha_0} \mathbb{E} \int_0^t |G_\alpha(Z_s^x) - G_{\alpha_0}(Z_s^x)|^2 \, ds = 0,$$

hence, $\exp \rho^\alpha(Z^x) \to \exp \rho^{\alpha_0}(Z^x)$ $\mathbb{P}$-a.s. In order to prove uniform integrability of $(\exp \rho^\alpha(Z^x))$, $\alpha \in \mathcal{A}$, it is enough to show

(8.5) $$\sup_{\alpha \in \mathcal{A}} \mathbb{E} \rho_\alpha(Z^x) \exp \rho_\alpha(Z^x) < \infty.$$

Setting $\widehat{W}_s := W_s - \int_0^s G_\alpha(Z_\tau^x) \, d\tau$, $s \in [0, t]$, we obtain, in virtue of the Girsanov theorem

(8.6)
$$\mathbb{E}\left( \int_0^t \langle G_\alpha(Z_s^x), dW_s \rangle - \tfrac{1}{2} \int_0^t |G_\alpha(Z_s^x)|^2 \, ds \right) \exp \rho_\alpha(Z^x)$$
$$= \mathbb{E}\left( \int_0^t \langle G_\alpha(Z_s^x), d\widehat{W}_s \rangle + \tfrac{1}{2} \int_0^t |G_\alpha(Z_s^x)|^2 \, ds \right) \exp \rho_\alpha(Z^x)$$
$$= \mathbb{E} \tfrac{1}{2} \int_0^t |G_\alpha(X_s^\alpha)|^2 \, ds \leq K^2 + K^2 \mathbb{E} \int_0^t |X_s^\alpha|_E^{2m} \, ds \leq N,$$

where $N < \infty$ is a constant independent of $\alpha \in \mathcal{A}$, which may be easily seen similarly as in (6.2) (cf. also Proposition 2.1 of [14]). Thus, $(\exp \rho^\alpha(Z^x))$ are uniformly integrable, which concludes the proof of (8.4). Now we have

(8.7)
$$\|\mu_\alpha - \mu_{\alpha_0}\|_{\text{var}} \leq \|P^\alpha(t, x, \cdot) - \mu_\alpha\|_{\text{var}}$$
$$+ \|P^{\alpha_0}(t, x, \cdot) - P^{\alpha_0}(t, x, \cdot)\|_{\text{var}}$$
$$+ \|P^{\alpha_0}(t, x, \cdot) - \mu^{\alpha_0}\|_{\text{var}}$$

and by Theorem 6.3,

$$\lim_{t \to \infty} \sup_{\alpha \in \mathcal{A}} \|P^\alpha(t, x, \cdot) - \mu_\alpha\|_{\text{var}} = 0,$$

which together with (8.4) yields (8.3). $\square$

## 9. Examples.

EXAMPLE 9.1 (Finite-dimensional equation). In the finite-dimensional case $E = H = \mathbb{R}^d$, the condition (2.7) is satisfied, even if the covariance matrix $Q$ is degenerate, which may be shown by generalizing a well-known Seidman's result [38] (cf. [25], Theorem 5.25). Obviously, $Z \in C([0, T], \mathbb{R}^d)$ for each $T > 0$, so the only assumptions that are needed in Theorem 6.3 ($V$-uniform ergodicity) and, if $P_t = P_t^*$, in Corollary 7.4 (spectral gap), are the strong Feller property for the linear equation (2.4) (which is true if and only if the matrix $Q_1$ is positive and is implied by positivity of the matrix $Q$), Hypothesis 2.4 and the the ultimate boundedness of solutions to (6.11). In order to apply Theorem 6.4 (uniform ergodicity) and Theorem 7.2 [spectral gap in $L^p(E,\mu^*)$], we have to assume the stronger growth condition (6.3).

As a specific example, we consider a nonlinear stochastic oscillator equation

$$\ddot{y} = f(y,\dot{y}) + \sigma \dot{w}_t, \qquad y(0) = x_1, \qquad \dot{y}(0) = x_2, \tag{9.1}$$

in $\mathbb{R}^d$. We assume that $f:\mathbb{R}^d \times \mathbb{R}^d \to \mathbb{R}^d$ is a locally Lipschitz function, $x_1, x_2 \in \mathbb{R}^d$, $\sigma \in \mathcal{L}(\mathbb{R}^d)$ is a regular matrix, and $(w_t)$ is a standard Wiener process in $\mathbb{R}^d$. Equation (9.1) may be rewritten in the form (2.1) with $X_t = (y(t), \dot{y}(t)) \in \mathbb{R}^{2d} = E = H$,

$$A = \begin{pmatrix} 0 & I \\ 0 & 0 \end{pmatrix}, \qquad F(x) = \begin{pmatrix} 0 \\ f(x) \end{pmatrix}, \qquad x \in \mathbb{R}^{2d},$$

$$Q^{1/2} = \begin{pmatrix} 0 & 0 \\ 0 & (\sigma\sigma^*)^{1/2} \end{pmatrix}.$$

According to the Kalman rank condition (see, e.g., [27]), the matrix $Q_t$ is invertible for each $t > 0$, so the equation with $F = 0$ is strongly Feller. Suppose that $f$ has at most polynomial growth, that is,

$$|f(x)|_{\mathbb{R}^d} \leq K(1 + |x|_{\mathbb{R}^{2d}}^m), \qquad x \in \mathbb{R}^{2d},$$

for some $K, m < \infty$. Then by Theorem 6.3, the solution is $V$-uniformly ergodic (Theorem 6.3), provided the ultimate boundedness condition (6.11) holds true. For example, we may take $d = 1$, $f(x) = -\alpha_2 x_2 - \alpha_1 x_1$ for $x = (x_1, x_2) \in \mathbb{R}^{2d}$, with some $\alpha_1, \alpha_2 > 0$ (a damped linear oscillator). Then (6.11) holds with constants which may be easily expressed in terms of $\alpha_1, \alpha_2$ and $\sigma$ and $V$-uniform ergodicity holds true. Similar results for a more general version of equation (9.1) can be found in [30].

Note that the semigroup $(P_t)$ is not symmetric in this case and Corollary 7.4 (on the spectral gap) is not applicable. Indeed, it follows from the results in [7, 8] that the spectral gap is zero in the present case. This example also shows that the assumption of symmetry of $P_t$ in Corollary 7.4 may not be removed.

EXAMPLE 9.2 (Stochastic reaction–diffusion equation with the cylindrical noise). Consider the system

$$\frac{\partial u}{\partial t} = Lu + f(u) + \eta,$$

$$u(0,\xi) = x(\xi), \tag{9.2}$$

$$\frac{\partial u}{\partial \xi}(t,0) = \frac{\partial u}{\partial \xi}(t,1) = 0, \qquad (t,\xi) \in \mathbb{R}_+ \times (0,1),$$

where $L$ is a uniformly elliptic operator

$$(9.3) \quad [L\varphi](\xi) = \left(\frac{\partial}{\partial \xi} a(\xi) \frac{\partial}{\partial \xi} \varphi\right)(\xi) + b(\xi) \frac{\partial}{\partial \xi} \varphi(\xi) + c(\xi) \varphi(\xi), \qquad \xi \in (0,1),$$

with $a, b, c \in C^1([0,1]), a(\xi) \geq a_0 > 0, \xi \in (0,1), f : \mathbb{R} \to \mathbb{R}$ is a locally Lipschitz mapping, and $\eta = \eta(t, \xi)$ is a nondegenerate noise. Let us note that the $C^1$ regularity of the coefficients is made for simplicity only and may be easily relaxed. The system (9.2) is rewritten in the form (2.1), with the coefficients defined in an obvious way on the spaces $H = L^2(0,1)$, $E = C([0,1])$,

$$A = L, \qquad \text{Dom}(A) = \left\{ \varphi \in H^2(0,1), \frac{\partial \varphi}{\partial \xi}(0) = \frac{\partial \varphi}{\partial \xi}(1) = 0 \right\}, \qquad Q \in \mathcal{L}(H),$$

where $Q$ is supposed to be boundedly invertible on $H$, and $F : E \to E, F(x(\xi)) = f(x(\xi)), \xi \in (0,1), x \in E$. It is well known (see, e.g., [9], A5.2) that $A$ generates a strongly continuous semigroup on $H$ and Hypothesis 2.2 is satisfied. With no loss of generality (replacing, if necessary, $A$ and $F$ by $A - \omega I$ and $F + \omega I$, respectively, with $\omega$ sufficiently large), we may assume that $\langle \tilde{A}x, x^* \rangle_{E, E^*} \leq 0$ for each $x \in \text{Dom}(\tilde{A}), x \in \partial \|x\|$ (recall that $\tilde{A}$ denotes the part of $A$ on $E$), then (6.1) is satisfied (see, e.g., [14], Example 3.1). It is well known that the corresponding Ornstein–Uhlenbeck process is strongly Feller and for a certain $c > 0$,

$$\|Q_t^{-1/2} S_t\| \leq \frac{c}{\sqrt{t}}, \qquad t \in (0,1),$$

(cf. [9]), hence, (3.12) holds. Moreover, standard estimates on the Green function of the problem [4] yield $\|S_t\|_{\text{HS}} \leq \text{const.} t^{-1/4}$, which implies (3.11) with $0 < \alpha < \frac{1}{2}$. Therefore, Hypotheses 2.1 and 2.3 hold and it remains only to specify the growth conditions on $f$. We assume that

$$(9.4) \qquad |f(\xi)| \leq k(1 + |\xi|^m), \qquad \xi \in \mathbb{R},$$

and

$$(9.5) \qquad f(\xi + \eta) \operatorname{sign} \xi \leq -k_1 |\xi| + k_2 |\eta|^s + k_3, \qquad \xi, \eta \in \mathbb{R},$$

for some constants $k, m, k_1, k_2, k_3$ and $s$. Now it is easy to check that Hypothesis 2.4 is satisfied and we may apply Theorem 6.3 to get $V$-uniform ergodicity with the rate which is specified there. Also, in the case the Markov semigroup $P_t$ is symmetric in $L^2(E, \mu^*)$ (e.g., if $Q = I$), we may apply Corollary 7.4 to obtain a lower bound for spectral gap. If the condition (9.5) is strengthened to

$$(9.6) \qquad f(\xi + \eta) \operatorname{sign} \xi \leq -k_1 |\xi|^{1+\varepsilon} + k_2 |\eta|^s + k_3, \qquad \xi, \eta \in \mathbb{R},$$

where $\varepsilon > 0$, then (6.3) holds as well and we may apply Theorem 6.4 on uniform exponential ergodicity and Theorem 7.2 on the spectral gap in $L^p(E, \mu^*), p \in [1, \infty)$. For example, if $f$ is a true polynomial, we have obtained the following result:

COROLLARY 9.3. *In Example* 9.2, *assume that* $\Psi$ *is a set of polynomials of the form*

$$f(\xi) = -a_{2n+1}\xi^{2n+1} + \sum_{i=0}^{2n} a_i \xi^i,$$

*where* $a_i$, $i = 1, \ldots, 2n$, *are in a given bounded set in* $\mathbb{R}^{2n}$, $a_{2n+1} \geq \bar{a}$, *for a given* $\bar{a} > 0$ *and* $n \geq 0$. *Then the* $V$-*uniform ergodicity (and if* $n > 0$, *uniform exponential ergodicity) holds with constants in* (6.24), (6.26) *and* (6.30) *uniform with respect to* $f \in \Psi$. *Also, if* $n > 0$, *then there is a positive lower bound on the spectral gap for* $(P_t)$ *in* $L^p(E, \mu^*), p \in (1, \infty)$, *uniform in* $f \in \Psi$. *Finally, if* $b = 0$ *and* $Q = I$, *then this bound holds also for* $p = 1$.

EXAMPLE 9.4 (The case of Lipschitz drift). Consider (9.2) with the same differential operator $L$ and initial and boundary conditions in the case when the noise may degenerate, for simplicity, suppose that $c \leq c_0 < 0$. For $\sigma \geq 0$, let $H_\sigma$ denote the domain $\text{dom}((-A)^\sigma)$ equipped with the graph norm $|y|_\sigma := |(-A)^\sigma y|$. As well known, for $\sigma \in (0, \frac{1}{2})$, the norm $|\cdot|_\sigma$ is equivalent with the norm of Sobolev–Slobodetskii space $H^{2\sigma}(0, 1)$,

$$|y|^2_{H^{2\sigma}} := |y|^2_{L^2(0,1)} + \int_0^1 \int_0^1 \frac{|y(\xi) - y(\eta)|^2}{(\xi - \eta)^{1+2\sigma}} \, d\xi \, d\eta, \qquad y \in H_\sigma.$$

Assume that $f : \mathbb{R} \to \mathbb{R}$ is Lipschitz continuous. It is easy to check that $F : H_\sigma \to H_\sigma$ is continuous and satisfies the growth condition (2.11) if, for some $\bar{k}_1, \bar{k}_2 > 0$, we have

$$(f(\xi) - f(\eta))\operatorname{sign}(\xi - \eta) \leq \bar{k}_1 |\xi - \eta| + \bar{k}_2, \qquad \xi, \eta \in \mathbb{R}.$$

Assume that $Q = (-A)^{-2\Delta}$ for some $\Delta \geq 0$. Then setting $E = H_\sigma$, we obtain

$$|G(x)| = |Q^{-1/2} F(x)| \leq \|Q^{-1/2}\|_{\mathcal{L}(H_\sigma, H)} |F(x)_\sigma \leq \widehat{K}(1 + |x|_\sigma), \qquad x \in H_\sigma,$$

for a suitable $\widehat{K} < \infty$, provided $\Delta \leq \sigma$. Moreover, the mapping $G : H_\sigma \to H$ is continuous since $F$ is continuous in $E$. In view of Remark 2.5, the results of the paper can be applied to this case. Also, condition (3.11) is satisfied with $\alpha < \frac{1}{2}$, as shown in Example 9.2, and Hypothesis 2.2(b) holds true for $E = H_\sigma$, provided $\sum \alpha_i^{2\sigma - 2\Delta - 1} < \infty$, where $(\alpha_i)$ are the eigenvalues of the operator $(-A)$ (cf. [9]), which is true (taking into account that $\alpha_i \sim i^2$) if $\Delta > \sigma - \frac{1}{4}$. The remaining condition (2.7) is always satisfied because $Q_t^{-1/2} S_t Q^{1/2} = \overline{Q^{1/2} S_t Q_t^{-1/2}} = \sqrt{2}(-A)^{1/2}(I - e^{tA})^{-1/2} e^{tA}$ and by the previous Example 9.2, we have that $\|H_t\|_{\text{HS}} \leq \text{const.} t^{-3/4}$. Summarizing, assume that

$$\sigma - \tfrac{1}{4} < \Delta \leq \sigma < \tfrac{1}{2}$$

holds, which may be achieved by a suitable choice of $\sigma \in (0, \frac{1}{2})$ for $\Delta \in [0, \frac{1}{2})$. Then $V$-uniform ergodicity follows from Theorem 6.3. Moreover, if $b = 0$, then the existence of the spectral gap in $L^2(E, \mu^*)$ with $E = H_\sigma$ follows from Corollary 7.4. In both cases the estimates on the rate of convergence are specified in Theorem 6.3 and Corollary 7.4, respectively.

## APPENDIX

For the reader's convenience, we collect here some basic facts about measurable linear mapping that are used in the paper. Most of them are well known.

Let $H$ be a real separable Hilbert space and let $\mu = N(0, C)$ be a centered Gaussian measure on $H$ with the covariance operator $C$ such that $\overline{\mathrm{im}(C)} = H$. The space $H_C = \mathrm{im}(C^{1/2})$ endowed with the norm $|x|_C = |C^{-1/2}x|$ can be identified as the reproducing kernel Hilbert space of the measure $\mu$. In the sequel we will denote by $\{e_n : n \geq 1\}$ the eigenbasis of $C$ and by $\{c_n : n \geq 1\}$ the corresponding set of eigenvalues:

$$Ce_n = c_n e_n, \qquad n \geq 1.$$

For any $h \in H$, we define

$$\phi_n(x) = \sum_{k=1}^n \frac{1}{\sqrt{c_k}} \langle h, e_k \rangle \langle x, e_k \rangle, \qquad x \in H.$$

LEMMA A.1. *The sequence $(\phi_n)$ converges in $L^2(H, \mu)$ to a limit $\phi$ and*

$$\int_H |\phi(x)|^2 \mu(dx) = |h|^2.$$

*Moreover, there exists a measurable linear space $\mathcal{M}_h \subset H$, such that $\mu(\mathcal{M}_h) = 1$, $\phi$ is linear on $\mathcal{M}_h$ and*

(A.1) $$\phi(x) = \lim_{n \to \infty} \phi_n(x), \qquad x \in \mathcal{M}_h.$$

*We will use the notation $\phi(x) = \langle h, C^{-1/2} x \rangle$.*

Let $H_1$ be another real separable Hilbert space and let $T : H \to H_1$ be a bounded operator. The Hilbert–Schmidt norm of $T$ will be denoted by $\|T\|_{\mathrm{HS}}$. Let

$$\tilde{T}_n x = \sum_{k=1}^n \frac{1}{\sqrt{c_k}} \langle x, e_k \rangle T e_k, \qquad x \in H.$$

LEMMA A.2. *Let $T: H \to H_1$ be a Hilbert–Schmidt operator. Then the sequence $(\tilde{T}_n)$ converges in $L^2(H, \mu; H_1)$ to a limit $\tilde{T}$ and*

$$\int_H |\tilde{T}(x)|_{H_1}^2 \mu(dx) = \|T\|_{\mathrm{HS}}^2.$$

*Moreover, there exists a measurable linear space $\mathcal{M}_T \subset H$, such that $\mu(\mathcal{M}_T) = 1$, $\tilde{T}$ is linear on $\mathcal{M}_T$ and*

(A.2) $$\tilde{T}(x) = \lim_{n \to \infty} \tilde{T}_n x, \qquad x \in \mathcal{M}_T.$$

*We will use the notation $TC^{-1/2}x = \tilde{T}(x)$.*

LEMMA A.3. *Let $K(t,s): H \to H$ be an operator-valued, strongly measurable function, such that, for each $a \in (0,1)$,*

(A.3) $$\int_0^1 \int_0^1 \|K(t,s)\|_{\mathrm{HS}} \, ds \, dt + \int_0^a \int_0^a \|K(t,s)\|_{\mathrm{HS}}^2 \, ds \, dt < \infty.$$

*Then the mapping $(t,s,y) \to K(t,s)C^{-1/2}y$ is measurable, and there exists a measurable linear space $\mathcal{M} \subset H$ of full measure, such that, for each $y \in \mathcal{M}$,*

$$\int_0^1 |K(t,s)C^{-1/2}y| \, ds < \infty, \qquad t\text{-a.e.}$$

PROOF. Let

$$\mathcal{K}x = K(t,s)x, \qquad x \in H.$$

By assumption, the operator $\mathcal{K}: H \to L^2((0,a) \times (0,a); H)$ is Hilbert–Schmidt for any $a < 1$ and thereby, by Lemma A.2, there exists the space $\mathcal{M}_a$ of full measure such that

$$\mathcal{K}C^{-1/2}y = \sum_{k=1}^\infty \frac{1}{\sqrt{c_k}} \langle y, e_k \rangle \mathcal{K}e_k,$$

where the convergence holds in mean-square and for each $y \in \mathcal{M}_a$, in $L^2((0,a) \times (0,a); H)$. Therefore, $\mathcal{K}C^{-1/2}$ is a measurable function of $(y,s,t)$ for $s,t \leq a$. Let $a_n \to a$ be increasing and let $\mathcal{M} = \bigcap_{n=1}^\infty \mathcal{M}_{a_n}$. Then $\mathcal{K}C^{-1/2}y$ is well defined, for all $s,t < 1$ is clearly measurable in $(y,s,t)$. Moreover, for each $y \in \mathcal{M}$,

$$I_n^2(y) = \left( \int_0^{a_n} \int_0^{a_n} |K(t,s)C^{-1/2}y| \, ds \, dt \right)^2$$
$$[3pt] \leq \int_0^{a_n} \int_0^{a_n} |K(t,s)C^{-1/2}y|^2 \, ds \, dt < \infty.$$

Since the sequence $I_n(y)$ is nondecreasing to a limit $I_\infty(y)$ for each $y \in \Omega$ and

$$\int_H I_\infty(y) \mu(dy) \leq \int_0^1 \int_0^1 \|K(t,s)\|_{\mathrm{HS}} \, ds \, dt < \infty,$$

the lemma follows. $\square$

School of Mathematics  
University of New South Wales  
Sydney 2052  
Australia  
E-mail: B.Goldys@unsw.edu.au

Mathematical Institute  
Academy of Sciences of Czech Republic  
Žitná 25  
11567 Prague 1  
Czech Republic  
E-mail: maslow@math.cas.cz